\documentclass[12pt, leqno]{amsart}

\usepackage{mathtools}
\mathtoolsset{showonlyrefs}
\numberwithin{equation}{section}

\usepackage{hyperref, graphicx}

\usepackage{amssymb}
\usepackage{amsthm}
\usepackage{amscd}
\usepackage{amsmath}
\usepackage{epic,eepic}
\usepackage{mathrsfs}
\usepackage{latexsym}
\usepackage{pifont}
\usepackage[all]{xy}
\usepackage {cancel}
\usepackage{verbatim}

\theoremstyle{plain}
\newtheorem{lemma}{Lemma}[section]
\newtheorem{prop}[lemma]{Proposition}
\newtheorem{thm}[lemma]{Theorem}
\newtheorem{cor}[lemma]{Corollary}
\newtheorem{intthm}{Theorem}

\theoremstyle{definition}

\newtheorem{rem}[lemma]{Remark}
\newtheorem{defi}[lemma]{Definition}
\newtheorem{exa}[lemma]{Example}

\newcommand{\bde}{\begin{defi}}
\newcommand{\ede}{\end{defi}\vspace{1mm}}
\newcommand{\ble}{\begin{lemma}}
\newcommand{\ele}{\end{lemma}}
\newcommand{\bpr}{\begin{prop}}
\newcommand{\epr}{\end{prop}}
\newcommand{\bt}{\begin{thm}}
\newcommand{\et}{\end{thm}}
\newcommand{\bco}{\begin{cor}}
\newcommand{\eco}{\end{cor}}
\newcommand{\bre}{\begin{rem}}
\newcommand{\ere}{\end{rem}}
\newcommand{\bex}{\begin{exa}}
\newcommand{\eex}{\end{exa}}
\newcommand{\bpf}{\begin{proof}}
\newcommand{\epf}{\end{proof}}

\newcommand{\mcB}{\mathcal{B}}

\newcommand{\mcD}{\mathcal{D}}
\newcommand{\mcE}{\mathcal{E}}
\newcommand{\mcF}{\mathcal{F}}
\newcommand{\mcG}{\mathcal{G}}
\newcommand{\mcH}{\mathcal{H}}

\newcommand{\mcL}{\mathcal{L}}
\newcommand{\mcM}{\mathcal{M}}
\newcommand{\mcN}{\mathcal{N}}
\newcommand{\mcO}{\mathcal{O}}
\newcommand{\mcP}{\mathcal{P}}
\newcommand{\mcQ}{\mathcal{Q}}
\newcommand{\mcR}{\mathcal{R}}
\newcommand{\mcS}{\mathcal{S}}
\newcommand{\mcT}{\mathcal{T}}
\newcommand{\mcU}{\mathcal{U}}

\newcommand{\mcW}{\mathcal{W}}

\newcommand{\mbG}{\mathbb{G}}

\newcommand{\mbP}{\mathbb{P}}
\newcommand{\mbQ}{\mathbb{Q}}
\newcommand{\mbR}{\mathbb{R}}

\newcommand{\mbZ}{\mathbb{Z}}

\newcommand{\mfa}{\mathfrak{a}}

\newcommand{\mfg}{\mathfrak{g}}

\newcommand{\msF}{\mathscr{F}}

\newcommand{\msN}{\mathscr{N}}

\newcommand{\mr}{\mathrm}
\newcommand{\N}{N}

\begin{document}

\title[The generic degree of the  generalized Verschiebung   in  rank two]{Explicit computation of the generic degree of \\ the generalized Verschiebung in rank two}
\author{Yuki Kondo}
\author{Yasuhiro Wakabayashi}

\address{\emph{Yuki Kondo}
 \newline
 \textnormal{Department of Mathematics, Graduate School of Science, The University of Osaka, Toyonaka, Osaka 560-0043, JAPAN.}
 \newline
 \textnormal{\texttt{u539693e@ecs.osaka-u.ac.jp}}}

 \address{\emph{Yasuhiro Wakabayashi}
 \newline
 \textnormal{Graduate School of Information Science and Technology, The University of Osaka, Suita, Osaka 565-0871, JAPAN.}
 \newline
  \textnormal{\texttt{wakabayashi@ist.osaka-u.ac.jp}}}

\date{}
\maketitle

\footnotetext{2020 {\it Mathematical Subject Classification}: Primary 14H60, Secondary 14G17;}
\footnotetext{Key words: positive characteristic, stable bundle, generalized Verschiebung, dormant oper}
\begin{abstract}
The purpose of this paper is to apply previous work on dormant opers to the study of the moduli space of stable bundles in positive characteristic. We affirmatively resolve the rank $2$ case of a conjecture proposed by the second author, which predicts a direct relationship between the number of higher-level dormant $\mathrm{PGL}_2$-opers and the generic degree of the generalized Verschiebung map for rank $2$ stable bundles induced by Frobenius pull-back. As a consequence, for a general smooth projective curve of genus $g>1$ in characteristic $p>2$, we obtain a procedure for explicitly determining these generic degrees by counting certain combinatorial objects and, in particular, compute $\mathrm{deg}(\mr{Ver}^2_1)$ as a polynomial in $p$ for $3\leq g\leq 10$.

\end{abstract}
\tableofcontents

\section{Introduction} \label{S0}

\subsection{Counting problem for  the generalized Verschiebung map} \label{SS100}
Let $k$ be an  algebraically closed field of characteristic $p >2$  and $X$  a  smooth   projective curve over $k$ of genus $g>1$.
For a positive integer $\N$,
 we denote by $X^{(\N)}$ the $\N$-th Frobenius twist of $X$ over $k$, i.e., 
 the base-change of $X$ along the $\N$-th iterate of the Frobenius automorphism of $\mr{Spec}(k)$.

 Pulling back  stable bundles on $X^{(\N)}$ via the $\N$-th relative  Frobenius morphism $F^{(\N)}_{X/k} : X \rightarrow X^{(\N)}$ gives rise to a rational map
\begin{align} \label{e1}
\mr{Ver}_{\N}^n : SU_{X^{(\N)}}^n \dashrightarrow SU_{X}^n 
\end{align}
 between the moduli spaces $SU_{X^{(\N)}}^n$ and  $SU_{X}^n$
   of rank $n >1$ stable bundles  with trivial determinant on $X^{(\N)}$ and $X$, respectively.
We  refer to this rational map as the {\bf ($\N$-th) generalized Verschiebung map}.
It captures the dynamics  of stable bundles under Frobenius pull-back and 
has been studied from various perspectives in, e.g., 
   ~\cite{BrKa}, ~\cite{dSa}, ~\cite{Gie1}, ~\cite{JoPa}, ~\cite{JRXY}, ~\cite{Las}, 
   \   ~\cite{LasPa1}, ~\cite{LasPa2},  
  ~\cite{Li}, and  ~\cite{Oss2}.
See also ~\cite{DuMe}, ~\cite{LanPa}, ~\cite{Li2}, ~\cite{Li3}, ~\cite{Li4}, ~\cite{Oss1},  and ~\cite{Yan}
 for studies on the geometric structure of  $\mr{Ver}^n_1$ itself.
 
  It is known that $\mr{Ver}^n_\N$ is  generically finite, which allows us to define its  generic degree $\mr{deg} (\mr{Ver}_\N^n)$, a quantity that reflects the  complexity of this rational map.
  Although the Frobenius action on the cohomology ring of the moduli stack of vector bundles has been investigated in ~\cite{CaNe} and  ~\cite{NeSt}, not much is known about  $\mr{deg} (\mr{Ver}_\N^n)$.
 It is therefore natural to ask for an explicit formula for this degree.
 Answering this question will lead to a deeper understanding of the nonabelian Hodge correspondence in positive characteristic, which arises from generalizations of  the transformations between vector bundles induced by  Frobenius pull-back (cf., e.g., ~\cite{GLQ}, ~\cite{OgVo}, ~\cite{Ohk}, ~\cite{Sch}).

Note that $\mr{Ver}^n_\N$ factors as a composite of $\N$ rational maps obtained as Frobenius twists of $\mr{Ver}^n_1$:
\begin{align} \label{Eq304}
SU_{X^{(\N)}}^n \dashrightarrow SU_{X^{(\N-1)}}^n  \dashrightarrow  \cdots  \dashrightarrow SU_{X^{(1)}}^n  \dashrightarrow SU_{X}^n.
\end{align}
This implies  the identity  $\mr{deg}(\mr{Ver}_\N^n) = \mr{deg}(\mr{Ver}^n_1)^{\N}$, which reduces the problem to the case $\N =1$.

Prior to this work, an explicit formula valid for all odd characteristics was known in the case $(n,g)=(2,2)$.
For a genus-$2$ curve $X$,  
the rational map $\mr{Ver}^2_1$ is relatively more accessible for analysis.
In fact, the moduli space $SU_{X}^2$ (as well as $SU_{X^{(1)}}^2$)  admits a natural compactification by semistable bundles, which  is canonically isomorphic to $\mbP^3$,  and its  boundary can be identified with the Kummer surface associated to $X$.
Moreover, the base points of the induced rational map between these compactifications correspond to Frobenius-unstable bundles. For each fixed theta characteristic on $X$, those bundles whose Frobenius pull-back has that theta characteristic as its destabilizing line subbundle are in bijection with {\bf dormant $\mr{PGL}_2$-opers} on $X$. 
Taking all $16$ theta characteristics into account, the enumeration of dormant opers established in ~\cite{Moc}, together with 
 ~\cite[Theorem 1.3]{Oss1} and  ~\cite[Corollary]{LanPa}, 
 yields the following formula for a general $X$:
 \begin{align} \label{Eq341}
 \mr{deg}(\mr{Ver}_{1}^{2}) = \frac{1}{3} \cdot (p^3+2p).
 \end{align}

For arbitrary genus  $g$,
the main result of ~\cite{HoWa} provides an upper bound for $\mr{deg}(\mr{Ver}^2_1)$, based on the idea of K. Joshi involving a correspondence  with certain Quot schemes (cf. ~\cite{Jos}, ~\cite{JoPa}).
However, this upper bound has been shown to be suboptimal (cf. the remark preceding ~\cite[Corollary 3.2.4]{HoWa}), and 
the overall structure of $\mr{Ver}^n_1$ remains mysterious (even in the case $n=2$) despite its fundamental significance.
In order to advance  the discussion beyond this point, 
we propose a new approach that connects the problem to the enumeration of certain combinatorial objects, via the theory of higher-level dormant opers developed in ~\cite{Wak4}.

\subsection{Comparison with dormant $\mr{PGL}_2$-opers of higher level} \label{SS101}

In the remainder of this introduction, 
we outline the contents of this paper.
A central object in our discussion is 
a {\bf dormant $\mr{PGL}_2^{(\N)}$-oper} (for $\N \in \mbZ_{> 0}$), which  is  a dormant $\mr{PGL}_2$-oper  enhanced by an action, with vanishing $p^\N$-curvature, 
 of  differential operators of level $\N -1$, introduced by P. Berthelot in the theory of arithmetic $\mcD$-modules (cf. ~\cite{Ber}).
Let  $\overline{\mcM}_g$ denote
 the moduli stack of genus-$g$ stable curves over $k$, and let
\begin{align} \label{Eq104}
\mcO p_{\N, g}^{^\mr{Zzz...}}
\end{align}
(cf. \eqref{Eq200}) denote
the moduli stack parametrizing pairs $(X, \msF^\spadesuit)$ consisting of a curve $X \in \overline{\mcM}_g$ and a dormant $\mr{PGL}_2^{(\N)}$-oper $\msF^\spadesuit$ on $X$.

As a generalization of ~\cite[Proposition 3.3]{JoXi} (see also ~\cite{JRXY}, ~\cite{JoPa}, ~\cite{Wak5}),
we show that
a higher-level version of Cartier descent yields (upon fixing a theta characteristic) a correspondence between dormant $\mr{PGL}_2^{(\N)}$-opers and certain rank $2$  stable bundles that are  called 
 {\bf maximally $F^{(\N)}$-destabilized} (cf.  Definition \ref{Def32}, Proposition \ref{Prop33}, and  \eqref{Eq55}).
Let $a$ be a $k$-rational point of $SU^2_{X^{(\N)}}$ corresponding to such a stable bundle, and
 let $b$ be the point in $\mcO p_{\N, g}^{^\mr{Zzz...}}$ determined by the associated dormant $\mr{PGL}_2^{(\N)}$-oper.
Then, one can verify  the existence of  a canonical isomorphism 
\begin{align}
\eta_\N : T_a SU^2_{X^{(\N)}} \xrightarrow{\sim} T_b \mcO p_{\N, g}^{^\mr{Zzz...}}
\end{align}
between the tangent spaces at $a$ and $b$, respectively (cf. Theorem \ref{Prop5}, (i)).
Moreover, for  each $\N' \in \mbZ_{> 0}$ with $\N' \leq \N$,  
the following square diagram commutes:
\begin{align} \label{Eq38rr}
\vcenter{\xymatrix@C=46pt@R=36pt{
 T_{a} SU^2_{X^{(\N)}} \ar[r]^-{\eta_\N}_-{\sim} \ar[d]_-{d\mr{Ver}^2_{\N \Rightarrow \N'} |_a} &  T_{b} \mcO p_{\N, g}^{^\mr{Zzz...}} \ar[d]^{d \Pi_{\N \Rightarrow \N'}|_{b}} \\
 T_{a'} SU^2_{X^{(\N')}}  \ar[r]^-{\sim}_-{\eta_{\N'}} &   T_{b'} \mcO p_{\N', g}^{^\mr{Zzz...}}
 }}
\end{align}
(cf. Theorem \ref{Prop5}, (ii), for the precise definitions of  the various notations).
In other words, 
the  first-order deformation theory of (the first $\N -1$ arrows of) the sequence \eqref{Eq304} is  compatible with that of the sequence
\begin{align} \label{Eq305}
\mcO p_{\N, g}^{^\mr{Zzz...}} \rightarrow \mcO p_{\N-1, g}^{^\mr{Zzz...}} \rightarrow \cdots \rightarrow \mcO p_{2, g}^{^\mr{Zzz...}} 
\rightarrow \mcO p_{1, g}^{^\mr{Zzz...}}
\end{align}
induced by reducing the level of dormant $\mr{PGL}_2$-opers.
This compatibility,  together with the generic \'{e}taleness of $\mcO p_{\N, g}^{^\mr{Zzz...}}/\overline{\mcM}_g$ established in ~\cite[Theorem C, (i)]{Wak4} and an extension argument based on very stability (cf. Section \ref{SS88}), leads to 
 the following statement, which forms  the first main result of this paper.

\begin{intthm}[cf. Corollary \ref{Cor88}] \label{ThA}
Let $\Pi_{\N}$ denote 
the natural projection $\mcO p_{\N, g}^{^\mr{Zzz...}} \rightarrow \overline{\mcM}_g$, which is 
 dominant  and generically finite  (cf. ~\cite[Theorems B and C, (i)]{Wak4}).
Then,  the following inequality holds:
\begin{align}  \label{Eq800}
\mr{deg} (\Pi_{\N}) \geq     \frac{\mr{deg}(\mr{Ver}_{1}^{2})^{\N -1} \cdot p^{g-1}}{2^{2g-1}} \cdot \sum_{\theta =1}^{p-1} \frac{1}{\sin^{2g-2} \big(\frac{\pi \cdot \theta}{p} \big)}.
\end{align}
Moreover, if $X$ is general (cf. Section \ref{SS1}), then this  inequality ``$\geq$" becomes an equality ``$=$".
In particular, the  resulting  formula gives an affirmative solution to   a strong form (for $n=2$)  of  the conjecture stated  in ~\cite[Conjecture  10.29]{Wak4} (cf. Remark \ref{Rem56}).
\end{intthm}

\subsection{Explicit computations via edge numberings on a trivalent graph} \label{SS102}

A key point of the preceding theorem is that the left-hand side of \eqref{Eq800} can be computed explicitly  by applying results from  ~\cite{Wak4}.
In fact, the decomposition of the degree  $\mr{deg}(\Pi_{\N})$ with respect to various clutching morphisms on moduli spaces of curves is  collectively described  through a    $2$d TQFT (= a $2$-dimensional topological quantum field theory).
This  $2$d TQFT provides  a combinatorial procedure for computing   $\mr{deg}(\Pi_{\N})$ by reducing 
the problem to the case where the underlying curve is totally degenerate.
In this degenerate case, each component of its normalization is isomorphic to a $3$-pointed projective line,  and dormant $\mr{PGL}_2^{(\N)}$-opers on such curves are determined entirely by combinatorial data associated with  their {\it radii} (cf. ~\cite[Chapter I, Definition 4.1]{Moc}).

Let us fix a finite connected  trivalent graph $\mbG$ of genus $g$ and  a totally degenerate curve  $X_\mbG$ whose dual graph coincides with $\mbG$ (cf. the discussion preceding Definition \ref{Def399}).
By a {\bf balanced $(p, \N)$-edge numbering} on $\mbG$, we mean  a collection of nonnegative integers $(a_e)_{e}$, indexed  by 
 the edges $e$ 
 of $\mbG$, such that for every  triple of edges $(e_1, e_2, e_3)$ (with multiplicity) incident to a common vertex, the triple $(a_{e_1}, a_{e_2}, a_{e_3})$ satisfies certain prescribed conditions, including specific types of triangle inequalities   (cf. Definition \ref{Def399} for the precise definition).
We denote the set of balanced  $(p, \N)$-edge numberings on $\mbG$  by 
\begin{align}
\mr{Ed}_{p, \N, \mbG}
\end{align}
(cf. \eqref{Eq111}).
This set is finite, and one can explicitly determine which combinations of nonnegative integers belong to it.

\vspace{8mm}
\hspace{0mm} 
 \includegraphics[width=16cm,bb=0 0 850 240,clip]{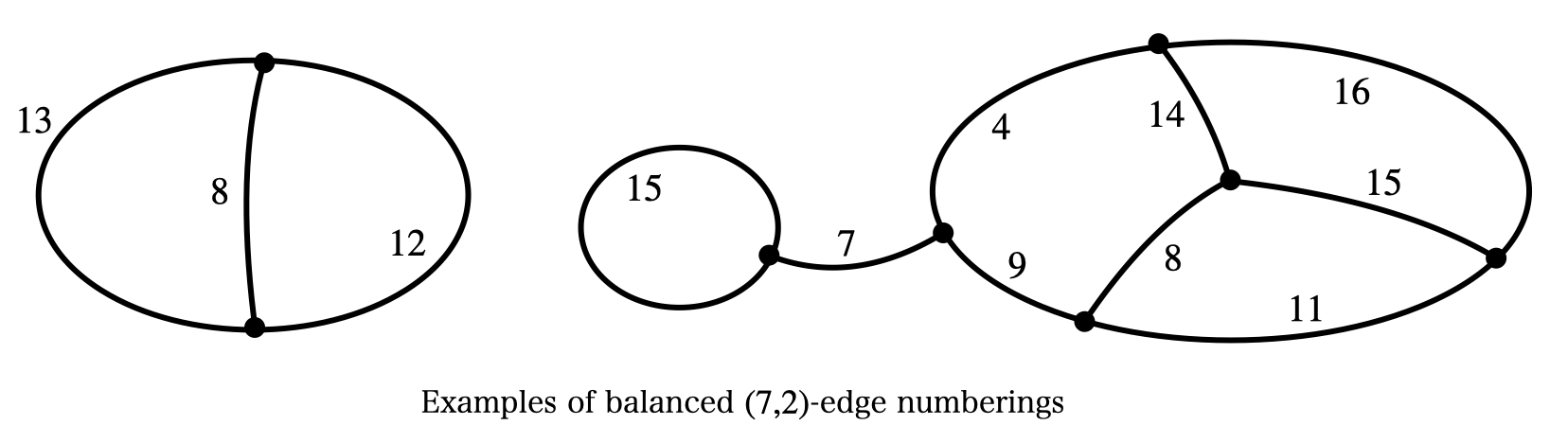}

\vspace{8mm}

According to  the theory of  diagonal reduction/lifting and  Gauss hypergeometric differential operators in characteristic $p^\N$, as  developed  in ~\cite{Wak4}, 
 dormant $\mr{PGL}_2^{(\N)}$-opers on $X_\mbG$
correspond bijectively to elements of $\mr{Ed}_{p, \N, \mbG}$.
  (This correspondence generalizes results for $\N =1$, discussed  in  ~\cite{LiOs}, ~\cite{Moc}, ~\cite{Wak2}, and ~\cite{Wak6}.)
In particular, Theorem \ref{ThA} implies the following result, which constitutes the second main result of this paper.

\begin{intthm}[cf. Theorems \ref{Th3} and \ref{Th56}] \label{ThE}
Suppose 
that $X$ is a general smooth projective curve classified by $\overline{\mcM}_g$.
Then, the following assertions hold:
\begin{itemize}
\item[(i)]
Let $\mbG$ be a finite connected trivalent graph of genus $g$.
Then, the following equality holds:
\begin{align} \label{Eq386}
 \mr{deg}(\mr{Ver}^2_1) = \frac{\sharp (\mr{Ed}_{p, 2, \mbG})}{\sharp (\mr{Ed}_{p, 1, \mbG})} \left(= \frac{\sharp (\mr{Ed}_{p, 2, \mbG})\cdot 2^{2g-1}}{p^{g-1}}  \cdot \left(\sum_{\theta =1}^{p-1} \frac{1}{\sin^{2g-2} \big(\frac{\pi \cdot \theta}{p} \big)}\right)^{-1} \right).
 \end{align}
 \item[(ii)]
 There exists a quasi-polynomial $Q (t)$ with coefficients in $\mbQ$, of degree $3g-3$ and period dividing $4$ (depending only on $g$),  that satisfies 
 \begin{align}
 \mr{deg} (\mr{Ver}_1^2) = Q (p).
 \end{align} 
 Moreover, the leading coefficient of $Q (t)$ coincides with $\frac{(-1)^g \cdot 2^{3g-4} \cdot B_{2g-2}}{(2g-2)!}$, where $B_{2g-2}$ denotes the $(2g-2)$-th Bernoulli number.
 \end{itemize}
\end{intthm}

As a consequence, the equality  \eqref{Eq386} allows for an explicit  computation of the generic degrees  by counting the number of  balanced $(p, \N)$-edge numberings for $\N =2$ (and $\N =1$).
Indeed,
 we recover not only the known formula 
 $\mr{deg} (\mr{Ver}^2_1) = \frac{p^3+ 2p}{3}$ for $g =2$ (as displayed  in \eqref{Eq341}), but also obtain further explicit expressions, as described in the following assertion.

\begin{intthm}[cf. Appendix \ref{AppPoly}] \label{ThG}
Suppose 
that $X$ is general.
Then, for each integer $g$ with $3 \leq g \leq 10$, the generic degree $\mr{deg}(\mr{Ver}^2_1)$ can be expressed as a polynomial in $p$.
The explicit formulas are listed in Appendix \ref{AppPoly}.
\end{intthm}

After the first version of this paper appeared, S. Zhang (cf. ~\cite{SZha}) used our combinatorial formula to prove that 
$\mr{deg}(\mr{Ver}_1^2)$ is a polynomial in $p$ for every genus $g\geq2$, and obtained an explicit formula.

\subsection{Notation and Conventions} \label{SS1}

Throughout this paper, we fix an integer $g > 1$, a prime number $p > 2$, and an algebraically closed field $k$ of characteristic $p$.
Also, fix a proper, smooth, and connected curve $X$ of genus $g$ over $k$.
Denote by 
 $\Omega_{X}$ the sheaf of K\"{a}hler differentials (i.e., $1$-forms) on $X$ over  $k$ and by $\mcT_X$ its dual, i.e.,  the tangent sheaf of $X/k$.

Let  $F_X$ denote  the absolute Frobenius endomorphism of $X$.
For a positive integer $\N$, the {\bf $\N$-th Frobenius twist} of $X$ over $k$ is, by definition, the base-change $X^{(\N)}$ of $X$ along the $\N$-th iterate of the absolute Frobenius endomorphism of $\mr{Spec}(k)$.
The morphism $F^{(\N)}_{X/k} : X \rightarrow X^{(\N)}$ induced naturally from $F^\N_X$ is called the {\bf $\N$-th relative Frobenius morphism} of
$X$ over $k$.
For convenience, we write $X^{(0)} := X$, $F^{(0)}_{X/k} := \mr{id}_X$, and $F_{X/k} := F_{X/k}^{(1)}$.

We shall write $\mfg := \mr{Lie} (\mr{SL}_2) \left(= \mr{Lie} (\mr{PGL}_2) \right)$, i.e., the Lie algebra of the special linear group $\mr{SL}_2$ over $k$ (identified with that of the projective linear group $\mr{PGL}_2$ because of the assumption $p >2$).
Given a rank $2$ vector bundle $\mcE$ on a curve $Y$ (such as $X$ or $X^{(\N)}$ for some $\N$),
we shall write $\mfg_\mcE$ for  the vector bundle associated to $\mcE$ via the adjoint representation $\mr{PGL}_2 \rightarrow \mr{Aut} (\mfg)$. 
Thus,  
$\mfg_\mcE$ is identified with   the sheaf $\mcE nd^0 (\mcE)$ of $\mcO_Y$-linear endomorphisms of $\mcE$ with vanishing trace.

Let   $\mcM_g$  denote the moduli stack of  proper, smooth, and  connected  curves of genus $g$  over $k$, and let  $\overline{\mcM}_g$ denote  its Deligne-Mumford  compactification, which parametrizes stable curves of genus $g$.
Since  $\mcM_g$   is 
an irreducible 
 Deligne-Mumford  stack over $k$ (cf. ~\cite[Section 5]{DeMu}), 
 it makes sense to speak of a ``{\it general}" curve, i.e., a curve corresponding to 
 a point of $\mcM_g$ that lies outside a fixed closed substack not equal to $\mcM_g$ itself.

For a rational map  $f : Y  \dashrightarrow Z$ between $k$-schemes,  with $Y$ integral and $Z$ separated over $k$,
we write  $\mr{Dom}(f)$ for  the  domain of definition of $f$, i.e., the largest open subscheme  of $Y$ on which $f$ is defined.


\vspace{10mm}
\section{Generalized Verschiebung maps} \label{S1}

In this section, we recall the definition of the generalized Verschiebung map between  moduli spaces of stable bundles and describe its  local structure   in terms of  cohomology groups.

\subsection{Moduli space of stable bundles} \label{SS2}

For a nonnegative integer $\N$, we denote by
\begin{align} \label{Eq30}
SU^2_{X^{(\N)}}
\end{align}
the coarse moduli space of rank $2$ stable bundles on $X^{(\N)}$ with trivial determinant. 
It is a smooth  irreducible variety  over $k$  of dimension $3g-3$ (cf. ~\cite[Lemma A]{MeSu}).

Given  another nonnegative integer $\N' \leq \N$,
the pull-back $\mcE \mapsto  F^{(\N - \N')*}_{X^{(\N')}/k}(\mcE)$ by $F^{(\N - \N')}_{X^{(\N')}/k}$ gives rise to  a rational map
\begin{align} \label{Eq31}
\mr{Ver}^2_{\N \Rightarrow N'} : SU^2_{X^{(\N)}}  \dashrightarrow  SU^2_{X^{(\N')}}
\end{align}
over $k$.
The domain  of definition $\mr{Dom}(\mr{Ver}^2_{\N \Rightarrow N'})$ classifies vector bundles $\mcE \in SU^2_{X^{(\N)}}$ whose pull-back 
 by  $F^{(\N-\N')}_{X^{(\N')}/k}$ remains stable.
 This description follows from the argument of ~\cite[Theorem A.6, (i)]{Oss1}, applied to the iterated relative Frobenius and restricted to the stable loci.
Since a vector bundle $\mcE'$ on $X^{(\N')}$ is stable whenever  its pull-back $F^{(\N' - \N'')*}_{X^{(\N'')}/k}(\mcE')$  is stable, for $0 \leq \N'' \leq \N'$,  
we obtain the following chain of inclusions: 
\begin{align} \label{Eq208}
\mr{Dom}(\mr{Ver}^2_{\N \Rightarrow 0}) \subseteq \mr{Dom}(\mr{Ver}^2_{\N \Rightarrow 1}) \subseteq \cdots \subseteq 
\mr{Dom}(\mr{Ver}^2_{\N \Rightarrow \N}) = SU^2_{X^{(\N)}}. 
\end{align}
For simplicity, we write $\mr{Ver}^2_{\N}  := \mr{Ver}^2_{\N \Rightarrow 0} $.

It  is well known that $\mr{Ver}^2_{\N \Rightarrow N'}$ is dominant and generically finite, so one may consider   its generic degree, denoted by $\mr{deg}(\mr{Ver}^2_{\N \Rightarrow N'})$.
To be precise, there exists a dense open subscheme $U_{\mr{ff}}$ of $SU_{X^{(\N')}}^2$ such that  $\mr{Ver}_{\N \Rightarrow \N'}^2$ restricts to a finite and faithfully flat morphism $(\mr{Ver}_{\N \Rightarrow \N'}^2)^{-1}(U_{\mr{ff}}) \rightarrow U_{\mr{ff}}$ (cf. ~\cite[Chapter II, Exercise 3.7]{Har}).
Moreover, if we assume that $X$ is  general in  $\mcM_g$, then  $\mr{Ver}^2_{\N \Rightarrow N'}$ is verified to be generically \'{e}tale (cf. ~\cite[Corollary 2.1.1]{MeSu}).
Under this assumption,  we can find a dense open subscheme $U_{\text{f\'{e}}}$ of $U_{\mr{ff}}$ such that  $\mr{Ver}_{\N \Rightarrow \N'}^2$ restricts to a finite and \'{e}tale morphism $(\mr{Ver}_{\N \Rightarrow \N'}^2)^{-1}(U_{\text{f\'{e}}}) \rightarrow U_{\text{f\'{e}}}$.

\subsection{Differential of the generalized Verschiebung map} \label{SS10}

Let $\N$ and $\N'$ be nonnegative integers with $\N' \leq \N$, and let
$a$ be a $k$-rational point of $\mr{Dom}(\mr{Ver}^2_{\N \Rightarrow \N'})$.
This point  determines  a stable bundle $\mcE$, whose pull-back $\mcE' := F^{(\N - \N')*}_{X^{(\N')}/k} (\mcE)$ by $F^{(\N - \N')}_{X^{(\N')}/k}$  remains stable.
Denote by $T_a SU^2_{X^{(\N)}}$ (resp., $T_{a'} SU^2_{X^{(\N')}}$) the tangent space of $SU^2_{X^{(\N)}}$ (resp., $SU^2_{X^{(\N')}}$) at $a$ (resp., $a' :=  \mr{Ver}_{\N \Rightarrow \N'}^2 (a)$).
From general  deformation theory of vector bundles,
we have 
a natural isomorphism of $k$-vector spaces
\begin{align} \label{Eq36}
\tau_{su} : T_a SU^2_{X^{(\N)}} \xrightarrow{\sim}  H^1 (X^{(\N)}, \mfg_\mcE) \ \left(\text{resp.,} \  
\tau'_{su} : T_{a'} SU^2_{X^{(\N')}}  \xrightarrow{\sim}  H^1 (X^{(\N')}, \mfg_{\mcE'})
\right).
\end{align}
Moreover, 
if $d\mr{Ver}^2_{\N \Rightarrow \N'} |_a$ denotes the $k$-linear morphism
$T_a SU^2_{X^{(\N)}} \rightarrow T_{a'} SU^2_{X^{(\N')}} $ obtained by differentiating $\mr{Ver}^2_{\N \Rightarrow \N'}$ at $a$,
then  the following square diagram commutes:
\begin{align} \label{Eq38}
\vcenter{\xymatrix@C=46pt@R=36pt{
 T_a SU^2_{X^{(\N)}} \ar[r]^-{\tau_{su}}_-{\sim} \ar[d]_-{d\mr{Ver}^2_{\N \Rightarrow \N'} |_a} &  H^1 (X^{(\N)}, \mfg_\mcE) \ar[d]^-{F^{(\N - \N')*}_{X^{(\N')}/k}} \\
 T_{a'} SU^2_{X^{(\N')}}  \ar[r]^-{\sim}_-{\tau'_{su}} &  H^1 (X^{(\N')}, \mfg_{\mcE'}).
 }}
\end{align}
Here,  the right-hand vertical arrow is  induced by pull-back along $F^{(\N - \N')}_{X^{(\N')}/k}$ and the canonical identification $F^{(\N - \N')*}_{X^{(\N')}/k}(\mfg_\mcE) \cong \mfg_{\mcE'}$.
Since both moduli spaces are smooth over $k$ of dimension $3g-3$,
the morphism  $\mr{Ver}^2_{\N \Rightarrow \N'}$ is \'{e}tale at the point $a$ if and only if this arrow  is an isomorphism.

\vspace{10mm}
\section{Dormant $\mr{PGL}_2$-opers of higher level} \label{S2}

This section  deals with  dormant  $\mr{PGL}_2$-opers of higher level, which were  introduced in ~\cite{Wak3} and  ~\cite{Wak4} as dormant $\mr{PGL}_2$-opers enhanced by actions of  differential operators of finite  level in the sense of  P. Berthelot. 
We also review a cohomological description of their deformation spaces, established  in ~\cite{Wak4} (cf. Proposition \ref{Prop18}).

\subsection{$\mcD$-modules of higher level} \label{SS4}

Fix a positive integer $\N$.
Following ~\cite{Ber}, 
the sheaf of differential operators $\mcD^{(\N -1)} := \mcD_{X/k}^{(\N -1)}$ on $X/k$  of level $\N -1$ is defined, where $\mr{Spec}(k)$ is equipped with the trivial $(\N -1)$-PD structure.

For each $j \in \mbZ_{\geq 0} \sqcup \{ \infty \}$, we write $\mcD^{(\N -1)}_{\leq j}$ for the subsheaf of $\mcD^{(\N -1)}$ consisting of differential operators of order  at most $j$, where $\mcD_{\leq \infty}^{(\N -1)} := \mcD^{(\N -1)}$.
In particular, we have $\mcD^{(\N -1)} = \bigcup_{j \in \mbZ_{\geq 0}} \mcD_{\leq j}^{(\N -1)}$.
Note that $\mcD_{\leq j}^{(\N -1)}$ carries  two distinct $\mcO_X$-module structures, i.e., one as given by left multiplication, where we denote this $\mcO_X$-module by ${^L}\mcD^{(\N -1)}_{\leq j}$, and the other given by right multiplication, where we denote this $\mcO_X$-module by ${^R}\mcD_{\leq j}^{(\N -1)}$.
For an $\mcO_X$-module $\mcE$, we equip the tensor product $\mcD_{\leq j}^{(\N -1)} \otimes \mcE  := {^R}\mcD_{\leq j}^{(\N -1)} \otimes \mcE$ with the $\mcO_X$-module structure given by left multiplication.

A {\bf (left) $\mcD^{(\N -1)}$-module structure} on $\mcE$ is a left $\mcD^{(\N -1)}$-action $\nabla : {^L} \mcD^{(\N -1)} \rightarrow \mcE nd_k (\mcE)$ on $\mcE$ extending its $\mcO_X$-module structure.
When it is necessary to emphasize the level ``$\N -1$", we write $\nabla^{(\N -1)}$ instead of $\nabla$.
A {\bf $\mcD^{(\N -1)}$-bundle} is defined as a pair $(\mcE, \nabla)$ consisting of a vector bundle $\mcE$ on $X$, and a $\mcD^{(\N -1)}$-module structure $\nabla$ on $\mcE$. 
For example, the pair $(\mcO_X, \nabla_{\mr{triv}}^{(\N -1)})$ of $\mcO_X$ and the natural $\mcD^{(\N -1)}$-module structure  $\nabla^{(\N -1)}_{\mr{triv}}$ on it forms a (line) $\mcD^{(\N -1)}$-bundle.
In the case of $\N =1$, giving a $\mcD^{(0)}$-module structure on an $\mcO_X$-module $\mcE$ amounts to giving a connection on $\mcE$, i.e., a $k$-linear morphism $\nabla : \mcE \rightarrow \Omega_X \otimes \mcE$ satisfying the Leibniz rule $\nabla (a \cdot v) = da \otimes v + a \cdot \nabla (v)$ for any $a \in \mcO_X$ and $v \in \mcE$.
(Since $X$ is of dimension $1$, any such connection is automatically flat, i.e., has vanishing curvature.)

Now, let   $(\mcE, \nabla)$ be
 a $\mcD^{(\N -1)}$-bundle. 
It is known that the  natural 
short exact sequence of $\mcO_X$-modules  $0 \rightarrow {^R}\mcD^{(\N -1)}_{\leq p^{\N}-1} \rightarrow {^R}\mcD^{(\N -1)}_{\leq p^{\N}} \rightarrow \mcT^{\otimes p^{\N}}_X\rightarrow 0$ has a canonical split injection 
$\psi_X : \mcT^{\otimes p^{\N}}_X \hookrightarrow \mcD^{(\N -1)}_{\leq p^{\N}}$
  (cf.  ~\cite[Section 2.5]{Wak4}).
The composite
  \begin{align} \label{Eq1030}
  \psi (\nabla) : \mcT^{\otimes p^{\N}}_X \xrightarrow{\psi_X}  \mcD^{(\N -1)}_{\leq p^\N} \xrightarrow{\mr{inclusion}}
    \mcD^{(\N -1)} \xrightarrow{\nabla} \mcE nd_{k} (\mcE),
\end{align}
 is called the {\bf $p^{\N}$-curvature} of $\nabla$ (cf. ~\cite{GLQ}, ~\cite[Definition 3.1.1]{LSQ}).
 When $\N =1$, the $p$-curvature is given by assigning $\partial^{\otimes p} \mapsto \nabla (\partial)^p - \nabla (\partial^{[p]})$ for each $\partial \in \mcT_X$, where $\partial^{[p]}$ denotes the local section of $\mcT_X$ corresponding to the derivation defined as  the $p$-th power iterate of $\partial$.
 
Denote by  $\mcS ol (\nabla)$  the subsheaf of $\mcE$ on which $\mcD_+^{(\N -1)}$ acts as zero via $\nabla$, where $\mcD_+^{(\N -1)}$ is the kernel of the canonical projection $\mcD^{(\N -1)} \twoheadrightarrow \mcO_X$.
 This sheaf can be regarded as an $\mcO_{X^{(\N)}}$-submodule of $F_{X/k*}^{(\N)} (\mcE)$ via  the underlying homeomorphism of $F_{X/k}^{(\N)}$, and sections in this sheaf are said to be  {\bf horizontal}.

For a vector bundle $\mcG$ on $X^{(\N)}$ (or more generally, an $\mcO_{X^{(\N)}}$-module),
 there exists a canonical
 $\mcD^{(\N -1)}$-module structure 
\begin{align} \label{E445}
\nabla_{\mcG, \mr{can}}^{(\N -1)} : {^L}\mcD^{(\N -1)} \rightarrow\mcE nd_k (F_{X/k}^{(\N)*}(\mcG))
\end{align}
on the pull-back $F_{X/k}^{(\N)*}(\mcG)$  of $\mcG$
along  $F_{X/k}^{(\N)}$ characterized  by the condition that every section of $(F^{(\N)}_{X/k})^{-1} (\mcG)$ is horizontal (cf. ~\cite[Corollaire 3.3.1]{Mon}).
The resulting assignments  $\mcG \mapsto (F_{X/k}^{(\N) *}(\mcG), \nabla_{\mcG, \mr{can}}^{(\N -1)})$
and $(\mcE, \nabla) \mapsto \mcS ol (\nabla)$ together 
 define an equivalence of categories
\begin{align} \label{Eq577}
\left(\begin{matrix} \text{the category of} \\ \text{vector bundles on $X^{(\N)}$} \end{matrix} \right) \xrightarrow{\sim}
\left(\begin{matrix} \text{the category of $\mcD^{(\N-1)}$-bundles} \\ \text{with vanishing $p^{\N}$-curvature} \end{matrix} \right)
\end{align}
 (cf. ~\cite[Corollary 3.2.4]{LSQ}), which  is compatible with standard operations such as direct sum ``$\oplus$", tensor product ``$\otimes$",  determinant ``$\mr{det}$'', and internal Hom  ``$\mcH om$". 
(The above discussion can be extended to the case where the underlying space is a family of smooth curves; see the proofs of the assertions in Section \ref{SS88}.)

Next, let $\N'$ be another positive integer with $\N' <  \N$, and $(\mcE, \nabla)$  a $\mcD^{(\N-1)}$-bundle with vanishing $p^\N$-curvature.
Note that there exists a natural morphism of sheaves of $k$-algebras  $\mcD^{(\N' -1)} \rightarrow \mcD^{(\N -1)}$ (cf. ~\cite[Section 2.2]{Ber}),
and that the composite
\begin{align} \label{Eq201}
\nabla^{(\N' -1)} : \mcD^{(\N' -1)} \rightarrow \mcD^{(\N -1)} \xrightarrow{\nabla} \mcE nd_k (\mcE)
\end{align}
specifies a $\mcD^{(\N' -1)}$-module structure on $\mcE$.
Since  $\nabla^{(\N' -1)}$ has  vanishing $p^{\N'}$-curvature,
it follows from  the equivalence of categories \eqref{Eq577} that the natural morphism
\begin{align} \label{Eq10}
F^{(\N -\N')*}_{X^{(\N')}/k} (\mcS ol (\nabla)) \rightarrow \mcS ol (\nabla^{(\N' -1)}) 
\end{align}
is an isomorphism between vector bundles on $X^{(\N')}$.

\subsection{Dormant $\mr{PGL}_2^{(\N)}$-opers and their moduli stack} \label{SS12}

Recall from ~\cite[Definitions 5.6, 5.7]{Wak4} that a {\bf dormant  $\mr{GL}_2^{(\N)}$-oper} on $X$ is 
a triple
\begin{align} \label{Eq1}
\msF^\heartsuit := (\mcF, \nabla, \mcL)
\end{align}
consisting of a $\mcD^{(\N -1)}$-bundle $(\mcF, \nabla)$ of rank $2$ and a line subbundle $\mcL$ of $\mcF$ satisfying the following conditions:
\begin{itemize}
\item
$\nabla$ has vanishing $p^\N$-curvature;
\item 
The composite 
\begin{align} \label{Eq40}
\mr{KS}_{\msF^\heartsuit} : \mcD_{\leq 1}^{(\N -1)} \otimes \mcL \xrightarrow{\mr{inclusion}} \mcD^{(\N -1)} \otimes \mcF \xrightarrow{\nabla} \mcF
\end{align}
is an isomorphism.
\end{itemize}
We refer to $\mcL$ as the {\bf Hodge bundle} of $\msF^\heartsuit$.
Also, a {\bf dormant $\mr{SL}_2^{(\N)}$-oper} is defined as  a dormant  $\mr{GL}_2^{(\N)}$-oper equipped with an isomorphism of line $\mcD^{(\N -1)}$-bundles $(\mr{det}(\mcF), \mr{det}(\nabla)) \xrightarrow{\sim} (\mcO_X, \nabla^{(\N -1)}_{\mr{triv}})$,
where $\mr{det}(\nabla)$ denotes the $\mcD^{(\N -1)}$-module structure on the determinant $\mr{det}(\mcF)$ induced naturally from $\nabla$.

The notion of an isomorphism between dormant $\mr{GL}_2^{(\N)}$-opers (resp., dormant $\mr{SL}_2^{(\N)}$-opers) is formulated in a natural manner.
Hence, one obtains the category of dormant $\mr{GL}_2^{(\N)}$-opers (resp., dormant $\mr{SL}_2^{(\N)}$-opers)  on $X$.

\begin{rem} \label{REm39}
Let $(\mcF, \nabla)$ be a $\mcD^{(\N -1)}$-bundle of rank $2$, and suppose that there exists a line subbundle $\mcL \subseteq \mcF$ such that  the triple $(\mcF, \nabla, \mcL)$ forms a dormant  $\mr{GL}_2^{(\N)}$-oper. 
The isomorphism  $\mr{KS}_{\msF^\heartsuit} : \mcD_{\leq 1}^{(\N -1)} \otimes \mcL \xrightarrow{\sim} \mcF$
  implies 
\begin{align}
\mr{deg} (\mcF/\mcL) & = \mr{deg} ((\mcD_{\leq 1}^{(\N -1)}\otimes \mcL)/(\mcD_{\leq 0}^{(\N -1)}\otimes \mcL)) \\
& 
= \mr{deg} (\mcT_X \otimes \mcL) \notag \\
& = 2-2g + \mr{deg} (\mcL) \notag \\
& < \mr{deg} (\mcL). \notag
\end{align}
This strict  inequality shows that such a line subbundle $\mcL$ is uniquely determined.
Therefore,   the groupoid of dormant $\mr{GL}_2^{(\N)}$-opers can be regarded as a full subgroupoid  of the groupoid of  $\mcD^{(\N -1)}$-bundles.
\end{rem}

Two dormant $\mr{GL}_2^{(\N)}$-opers $\msF^\heartsuit_i := (\mcF_i, \nabla_i, \mcL_i)$ ($i=1, 2$) on $X$ are said to be {\bf equivalent}
if there exists a line $\mcD^{(\N -1)}$-bundle $\msN := (\mcN, \nabla_\mcN)$ with vanishing $p^\N$-curvature  such that the twisted dormant $\mr{GL}_2^{(\N)}$-oper 
$(\msF_1^\heartsuit)_{\otimes \msN} := (\mcF_1 \otimes \mcN, \nabla_1 \otimes \nabla_\mcN, \mcL_1 \otimes \mcN)$ is isomorphic to  $\msF^\heartsuit_2$.
This binary relation forms  an equivalence relation on the set of dormant $\mr{GL}_2^{(\N)}$-opers on $X$, and we define a {\bf dormant $\mr{PGL}_2^{(\N)}$-oper} as an  equivalence class $\msF^\spadesuit := [\msF^\heartsuit]$ of a dormant $\mr{GL}_2^{(\N)}$-oper $\msF^\heartsuit$.

\begin{rem} \label{Rem9}
The notion of a dormant $\mr{PGL}_2^{(\N)}$-oper can be extended, within the framework of logarithmic geometry,  to the case where the underlying space is a family of  pointed stable curves (cf. ~\cite[Definitions 5.2,  5.3]{Wak4}).
This generalization  is essential for constructing  the  moduli stack over $\overline{\mcM}_g$ that will appear later  in \eqref{Eq200}, and for carrying out degeneration arguments that reduce various problems to the case of small genus (cf. Section \ref{SS5}).
However, for simplicity of our discussion, we omit the details of the generalized formulation   in this paper.
\end{rem}

Next, let us fix a theta characteristic  of $X$, i.e., a line bundle $\varTheta$ equipped with an isomorphism $\varTheta^{\otimes 2} \xrightarrow{\sim} \Omega_{X}$.
We shall write $\mcF_\varTheta := \mcD_{\leq 1}^{(\N -1)} \otimes \varTheta$, and regard
 $\varTheta$  as its  line subbundle via the natural identification $\varTheta = \mcD_{\leq 0}^{(\N -1)} \otimes \varTheta \left(\subseteq \mcF_\varTheta \right)$.
Then,  we have a sequence of isomorphisms
\begin{align} \label{Eq133}
\mr{det} (\mcF_\varTheta) \xrightarrow{\sim} \left(\mcF_\varTheta /\varTheta \right) \otimes \varTheta \xrightarrow{\sim}( \mcT_X \otimes \varTheta ) \otimes \varTheta \xrightarrow{\sim} \mcO_X,
\end{align}
where the final  arrow arises from the fixed isomorphism $\varTheta^{\otimes 2} \xrightarrow{\sim} \Omega_X \left(= \mcT_X^\vee \right)$.

\bde \label{Def2}
A {\bf dormant $(\mr{SL}^{(\N)}_2, \varTheta)$-oper}
  is a $\mcD^{(\N -1)}$-module structure $\nabla$ on $\mcF_\varTheta$ such that the triple 
$(\mcF_\varTheta, \nabla, \varTheta)$
forms a dormant $\mr{GL}_2^{(\N)}$-oper satisfying the equality $\mr{det}(\nabla) = \nabla_\mr{triv}^{(\N -1)}$ under the identification $\mr{det}(\mcF_\varTheta) = \mcO_X$ given by \eqref{Eq133}.
\ede

\begin{rem} \label{Rem44}
As explained in  ~\cite[Example 5.14]{Wak4}, the fixed theta characteristic $\varTheta$ defines a {\it $2^{(\N)}$-theta characteristic} $\vartheta_0$  in the sense of ~\cite[Definition 5.13]{Wak4}.
Then, the above definition coincides with the notion of a {\it dormant $(\mr{GL}_2^{(\N)}, \vartheta_0)$-oper} introduced in ~\cite[Definition 5.16, (i)]{Wak4}.
\end{rem}

Dormant $(\mr{SL}_2^{(\N)}, \varTheta)$-opers can be regarded as dormant $\mr{SL}_2^{(\N)}$-opers, so they form a full subcategory of dormant  $\mr{SL}_2^{(\N)}$-opers.
Moreover,  as shown in  ~\cite[Theorem 5.19]{Wak4},  taking equivalence classes yields  a bijection
\begin{align} \label{Eq55}
\left(\begin{matrix} \text{the set of isomorphism classes of} \\  
\text{dormant $(\mr{SL}_2^{(\N)}, \varTheta)$-opers  on $X$} \end{matrix} \right) \xrightarrow{\sim}
\left(\begin{matrix} \text{the set  of dormant} \\ \text{$\mr{PGL}_2^{(\N)}$-opers on $X$} \end{matrix} \right).
\end{align}

We denote by  
\begin{align} \label{Eq200}
\mcO p_{\N, g}^{^\mr{Zzz...}}
\end{align}
 the moduli stack classifying pairs $(X, \msF^\spadesuit)$ consisting of a curve $X$ classified by $\overline{\mcM}_g$ and a 
dormant $\mr{PGL}_2^{(\N)}$-oper $\msF^\spadesuit$ on $X$ (cf. Remark \ref{Rem9}).
According to ~\cite[Theorems B and C, (i)]{Wak4} (or ~\cite[Chapter II, Theorem 2.8]{Moc} in the case $\N =1$),
$\mcO p_{\N, g}^{^\mr{Zzz...}}$ is a nonempty smooth  Deligne-Mumford stack  of finite type over $k$ of dimension $3g-3$ and the natural projection 
\begin{align} \label{Eq250}
\Pi_{\N} : \mcO p_{\N, g}^{^\mr{Zzz...}} \rightarrow \overline{\mcM}_g
\end{align}
 is generically finite and generically \'{e}tale.
For  another positive integer  $\N' \leq  \N$,
we obtain a generically finite   morphism
\begin{align} \label{Eq25}
\Pi_{\N \Rightarrow \N'} : \mcO p_{\N, g}^{^\mr{Zzz...}} \rightarrow \mcO p_{\N', g}^{^\mr{Zzz...}}
\end{align}
over $\overline{\mcM}_g$ given  by reducing the levels of dormant $\mr{PGL}_2^{(\N)}$-opers (cf. \eqref{Eq201}).

\subsection{Infinitesimal deformations of a dormant $\mr{PGL}_2^{(\N)}$-oper} \label{SS34}

Let $\varTheta$ be as above.
For simplicity, we write $\mcF := \mcF_\varTheta$ and $\mcF^1 := \varTheta \left(\subseteq \mcF_\varTheta \right)$.
Consider the $3$-step decreasing filtration
\begin{align} \label{Eq220}
0= 
\mfg_\mcF^2 \subseteq \mfg_\mcF^1 \subseteq \mfg_\mcF^0 \subseteq \mfg_\mcF^{-1} = \mfg_\mcF
\end{align}
on the rank $3$ vector bundle $\mfg_\mcF$ defined in such a way  that $\mfg_\mcF^1$ (resp., $\mfg_\mcF^0$) consists of local endomorphisms  $h$ of $\mcF$ with $h (\mcF^1) = 0$ (resp., $h (\mcF^1) \subseteq \mcF^1$).
Then, the  associated graded pieces satisfy $\mfg_\mcF^j/\mfg_\mcF^{j+1} \cong \Omega_X^{\otimes j}$  for $j=-1, 0, 1$.

Now, we recall the Atiyah algebroid of the $\mr{SL}_2$-bundle corresponding to  $\mcF$.
Let us take an open subscheme $U$ of $X$ over which the inclusion $\varTheta \hookrightarrow \mcF$ can be identified with the inclusion into    the first  direct summand   $\mcO_U \hookrightarrow \mcO_U^{\oplus 2}$.
We choose this identification to be compatible with the determinant trivialization \eqref{Eq133}.
Under this identification, we define   $\widetilde{\mcT}_\mcF (U)$  (resp., $\widetilde{\mcT}_{\mcF, B} (U)$) to be the subset of  $\mr{End}_k (\mcF |_U) \left(= \mr{End}_k (\mcO_U^{\oplus 2})\right)$ consisting of endomorphisms of the form $h = D^{\oplus 2} + h'$, where
$D$ is  a derivation on $\mcO_U$ and $h' \in \Gamma (U, \mfg_\mcF)$ (resp., $h' \in \Gamma (U, \mfg_\mcF^0)$).
This subset does not depend on the choice of such a determinant-compatible identification.
We shall denote by
\begin{align} \label{Eq221}
\widetilde{\mcT}_{\mcF} \ \left(\text{resp.,} \  \widetilde{\mcT}_{\mcF, B} \right)
\end{align}
the sheaf associated to the assignment $U \mapsto \widetilde{\mcT}_\mcF (U)$ (resp., $U \mapsto \widetilde{\mcT}_{\mcF, B} (U)$).
If we consider  $\mcE nd_k (\mcF)$ as an $\mcO_X$-module via  left multiplication, then 
$\widetilde{\mcT}_{\mcF}$ (resp., $\widetilde{\mcT}_{\mcF, B}$)
defines an $\mcO_X$-submodule of this sheaf.
Moreover, it fits into the following  short exact sequence of $\mcO_X$-modules:
\begin{align} \label{Eq47}
0 \rightarrow \mfg_\mcF \rightarrow \widetilde{\mcT}_{\mcF} \xrightarrow{d_{\mcF}} \mcT_X \rightarrow 0 \ 
 \left(\text{resp.,} \  0 \rightarrow \mfg^0_\mcF \rightarrow \widetilde{\mcT}^{}_{\mcF, B}  \xrightarrow{d_{\mcF, B}} \mcT_X \rightarrow 0 \right).
\end{align}
See, e.g., ~\cite[Section 1.2.5]{Wak9} for a more  intrinsic construction of this sequence.

Next, let $\nabla$ be a dormant $(\mr{SL}_2^{(\N)}, \varTheta)$-oper on $X$.
The $\mcD^{(0)}$-module structure $\nabla^{(0)}$ induced from $\nabla$ (cf. \eqref{Eq201}) corresponds to a connection  on $\mcF$, hence yields  a split injection $\mcT_X \hookrightarrow \widetilde{\mcT}_{\mcF}$ of the non-resp'd short exact sequence  \eqref{Eq47}; we abuse the notation ``$\nabla^{(0)}$" for writing 
 this connection.

We shall denote by 
\begin{align} \label{Eq127}
\widetilde{\mcT}_{\mcF, \nabla}  \ \ \ \text{and} \ \ \ \mfg_{\mcF, \nabla}
\end{align}
(cf. ~\cite[Section 8.2]{Wak4}) the subsheaves of $\widetilde{\mcT}_\mcF  \left(\subseteq \mcE nd_k (\mcF) \right)$ and $\mfg_\mcF  \left(\subseteq \mcE nd_{\mcO_X} (\mcF) \right)$, respectively,  consisting of local sections $v$ satisfying the relation   
$[\nabla (D), v] - \nabla ([D, d_{\mcF} (v)]) = 0$  for any $D \in \mcD^{(\N -1)}$.
In particular, the equality $\mfg_{\mcF, \nabla} = \mfg_\mcF \cap \widetilde{\mcT}_{\mcF, \nabla}$ holds.
Also, we set
\begin{align}
\widetilde{\mcT}_{\mcF, B, \nabla} := \widetilde{\mcT}_{\mcF, B} \cap \widetilde{\mcT}_{\mcF, \nabla}.
\end{align}
These sheaves are regarded as $\mcO_{X^{(\N)}}$-modules via the underlying homeomorphism of $F^{(\N)}_{X/k}$.
 The surjection  $d_{\mcF, B} : \widetilde{\mcT}_{\mcF, B} \twoheadrightarrow \mcT_X$ induces 
 a morphism of $k$-vector spaces
 \begin{align} \label{Eq228ddd}
 D_{\nabla} : 
 H^1 (X, \widetilde{\mcT}_{\mcF, B, \nabla}) \rightarrow H^1 (X, \mcT_X).
 \end{align}
 According to ~\cite[Corollary 8.11, (ii)]{Wak4},  the following assertion holds.

 \bpr \label{Prop18}
 Denote by $\msF^\spadesuit$ the dormant $\mr{PGL}_2^{(\N)}$-oper corresponding to $\nabla$ via \eqref{Eq55} and by
 $b$ (resp., $\overline{b}$) the $k$-rational point of $\mcO p_{\N, g}^{^\mr{Zzz...}}$  (resp., $\overline{\mcM}_g$) classifying  $(X, \msF^\spadesuit)$ (resp., $X$).
 Also,  denote  by $T_b \mcO p_{\N, g}^{^\mr{Zzz...}}$ the tangent space of $\mcO p_{\N, g}^{^\mr{Zzz...}}$ at the point $b$.
 Then, there exists a canonical isomorphism of $k$-vector spaces
 \begin{align} \label{Eq202}
 \tau_{op} : T_b \mcO p_{\N, g}^{^\mr{Zzz...}} \xrightarrow{\sim} H^1 (X, \widetilde{\mcT}_{\mcF, B, \nabla}),
 \end{align}
 which makes the following square diagram commute:
 \begin{align} \label{Eq62}
\vcenter{\xymatrix@C=46pt@R=36pt{
T_b \mcO p_{\N, g}^{^\mr{Zzz...}} \ar[r]_-{\sim}^-{\tau_{op}} \ar[d]_-{d \Pi_{\N} |_b} &
H^1 (X, \widetilde{\mcT}_{\mcF, B, \nabla})
 \ar[d]^-{D_\nabla} \\
T_{\overline{b}} \overline{\mcM}_g \ar[r]^-{\sim}_-{\tau_m} & H^1 (X, \mcT_X),
 }}
\end{align}
where  $\tau_m$ denotes the classical  Kodaira-Spencer isomorphism at $\overline{b}$ and the left-hand vertical arrow $d \Pi_{\N} |_b$ is the differential of $\Pi_{\N}$ at $b$.
\epr

 Let us choose another positive  integer  $\N' \leq \N$.
 Note that the vector bundle $\mcF$ does not depend on the level ``$\N$".
 Then, the $\mcD^{(\N' -1)}$-module structure  $\nabla' := \nabla^{(\N' -1)}$ associated to $\nabla$ defines a dormant $(\mr{SL}^{(\N')}_2, \varTheta)$-oper.
 There exists a natural inclusion $\widetilde{\mcT}_{\mcF, B, \nabla} \hookrightarrow \widetilde{\mcT}_{\mcF, B, \nabla'}$, which  induces a morphism
 $\iota_{\N \Rightarrow \N'} : H^1 (X, \widetilde{\mcT}_{\mcF, B, \nabla}) \rightarrow H^1 (X, \widetilde{\mcT}_{\mcF, B, \nabla'})$.
 This morphism  fits into the following commutative diagram:
 \begin{align} \label{Ed63}
\vcenter{\xymatrix@C=46pt@R=36pt{
T_b \mcO p_{\N, g}^{^\mr{Zzz...}} \ar[r]_-{\sim}^-{\tau_{op}} \ar[d]_-{d \Pi_{\N \Rightarrow \N'}|_{b}} &
H^1 (X, \widetilde{\mcT}_{\mcF, B, \nabla})
\ar[d]^-{\iota_{\N \Rightarrow \N'}} \\
T_{b'} \mcO p_{\N', g}^{^\mr{Zzz...}}  \ar[r]^-{\sim}_-{\tau'_{op}} & 
H^1 (X, \widetilde{\mcT}_{\mcF, B, \nabla'}).
 }}
\end{align}
Here,  $b'$ denotes the image of $b$ under $\Pi_{\N \Rightarrow \N'}$,
$\tau'_{op}$ is the isomorphism ``$\tau_{op}$" corresponding to $\nabla'$, and 
 the left-hand vertical arrow $d \Pi_{\N \Rightarrow \N'}|_{b}$ denotes the differential of $\Pi_{\N \Rightarrow \N'}$ at $b$.

\vspace{10mm}
\section{Correspondence between stable bundles and dormant opers} \label{S2dd}

This section discusses  a correspondence between  dormant $\mr{SL}_2^{(\N)}$-opers and   
  certain vector bundles called   maximally $F^{(\N)}$-destabilized bundles (cf. Definition \ref{Def32}).
  This correspondence extends to a bijection between the respective first-order deformation spaces (cf. Theorem \ref{Prop5}).
 As a consequence of this result, combined with the generic \'{e}taleness of $\Pi_{\N}$,
 we will establish a connection  between the (generic) degrees of $\mr{Ver}^2_{1}$ and $\Pi_\N$ (cf. Corollary \ref{Cor88}).

\subsection{Maximally Frobenius-destabilized bundles} \label{SS8}

Let $\N$ be a positive integer and 
 $\mcE$  a vector bundle on $X$ of rank $2$.
The {\bf first Segre invariant} of $\mcE$ is defined as  
\begin{align}
s_1 (\mcE) := \mr{deg}(\mcE) - 2 \cdot  \mr{max} \left\{\mr{deg}(\mcL) \, | \, \text{$\mcL$ is a line subbundle of $\mcE$} \right\}.
\end{align}
In particular,
 $s_1 (\mcE) > 0$ (resp., $s_1 (\mcE) \geq 0$) if and only if $\mcE$ is stable (resp., semistable).
Moreover, if   $\mcE$ is indecomposable, then
$s_1 (\mcE) \geq  2-2g$.
Indeed, suppose (for contradiction) that $s_1 (\mcE) <  2-2g$.
Choose a line subbundle $\mcL$ of $\mcE$ with $s_1(\mcE) = \mr{deg}(\mcE) - 2 \cdot \mr{deg}(\mcL)$. Then, the assumption implies
\begin{align}
\mr{deg}(\mcL \otimes (\mcE/\mcL)^\vee) = 2 \cdot \mr{deg}(\mcL)   - \mr{deg}(\mcE)  = -s_1 (\mcE) > 2g-2.
\end{align}
It follows that
\begin{align}
\mr{Ext}^1 (\mcE/\mcL, \mcL) \cong H^1 (X, \mcL \otimes (\mcE/\mcL)^\vee) = 0,
\end{align}
so we have  $\mcE \cong \mcL \oplus (\mcE/\mcL)$, which is a contradiction.
With this in mind, we recall  the following definition (cf.  ~\cite{JoPa}, ~\cite{JRXY}, ~\cite{LanPa}, and ~\cite{Zha} for the case $\N =1$).

\bde[cf. ~\cite{Wak5}, Definition 33] \label{Def32}
Let $\mcE$ be a  stable  (and hence indecomposable) vector bundle of rank $2$ on $X^{(\N)}$.
We shall say that $\mcE$ is {\bf maximally $F^{(\N)}$-destabilized} if $F^{(\N -1)*}_{X^{(1)}/k}(\mcE)$ is semistable and 
the equality $s_1 (F^{(\N)*}_{X/k}(\mcE)) = 2-2g$ holds.
\ede

Let $\mcE$ be a (stable and) maximally $F^{(\N)}$-destabilized bundle of rank $2$ with trivial determinant.
By definition, there exists a line subbundle $\mcL \subseteq F^{(\N)*}_{X/k} (\mcE)$ of degree $g-1$, which is uniquely determined.
Note that $\mcL$ is {\it not} closed under 
the connection $\nabla$ on $F^{(\N)*}_{X/k} (\mcE)$ induced by 
the canonical $\mcD^{(\N -1)}$-module structure  $\nabla_{\mcE, \mr{can}}^{(\N -1)}$.
Indeed, if $\nabla (\mcL) \subseteq \Omega_X \otimes \mcL$,
then $\nabla$ restricts to a connection on $\mcL$, which has vanishing $p$-curvature.
Cartier descent  (i.e., the equivalence of categories \eqref{Eq577} for $\N =1$) gives a line subbundle $\mcL_1$ of  $F^{(\N -1)*}_{X^{(1)}/k} (\mcE)$ with $F^*_{X/k} (\mcL_1) \cong \mcL$.
Hence, $\mr{deg} (\mcL_1) = \frac{1}{p} \cdot \mr{deg}(\mcL) = \frac{g-1}{p} > 0$, which contradicts the semistability of $F^{(\N -1)*}_{X^{(1)}/k} (\mcE)$.

It follows  that the image of the composite
\begin{align} \label{hiwpss}
\phi : \mcD^{(0)}_{\leq 1} \otimes \mcL \xrightarrow{\mr{inclusion}} \mcD^{(0)} \otimes F^{(\N)*}_{X/k} (\mcE) \xrightarrow{\nabla} F^{(\N)*}_{X/k} (\mcE)
\end{align}
is locally free of  rank $2$.
Since $\phi(\mcD^{(\N -1)}_{\leq 0} \otimes \mcL ) = \mcL$, the morphism $\phi$ induces a nonzero morphism
\begin{align}
\bar{\phi} : \mcT_X \otimes \mcL \to F^{(\N)*}_{X/k} (\mcE) / \mcL .
\end{align}
Comparing the degrees of line bundles, we see that $\overline{\phi}$ is an isomorphism.
Thus,  $\phi$ is an isomorphism.

The determinant of the isomorphism $\phi$ yields
\begin{align}
 \mcT_X \otimes \mcL^{\otimes 2}  \xrightarrow{\sim} \mr{det} (\mcD^{(\N -1)}_{\leq 1} \otimes \mcL) \xrightarrow{\sim} \mr{det}(F^{(\N)*}_{X/k} (\mcE)) \xrightarrow{\sim}F^{(\N)*}_{X/k} (\mr{det}(\mcE)) \xrightarrow{\sim}F^{(\N)*}_{X/k}  (\mcO_{X^{(\N)}}) \xrightarrow{\sim} \mcO_X.
\end{align}
This implies   that $\mcL$ defines a theta characteristic of $X$.

Let $\nabla_\mcE$ denote the $\mcD^{(\N -1)}$-module structure on $\mcF_\mcL \left(:= \mcD_{\leq 1}^{(0)} \otimes \mcL \right)$ corresponding to $\nabla_{\mcE, \mr{can}}^{(\N -1)}$ via  \eqref{hiwpss}.
The determinant of $\nabla_\mcE$ is trivial because of the equivalence of categories  \eqref{Eq577} and $\mr{det}(\mcE) \cong \mcO_{X^{(\N)}}$.
Thus, the resulting triple 
\begin{align} \label{Eq209}
\msF^\heartsuit_\mcE := (\mcF_\mcL, \nabla_\mcE, \mcL)
\end{align}
 forms a dormant $\mr{SL}_2^{(\N)}$-oper, and moreover $\nabla_\mcE$ forms a dormant $(\mr{SL}_2^{(\N)}, \mcL)$-oper.

\begin{rem} \label{Rem33}
Let $\mcE$ be a stable bundle classified by a point $a \in SU^2_{X^{(\N)}} (k)$, and suppose that it is maximally $F^{(\N)}$-destabilized.
Then,  for each positive integer $\N' \leq \N$,
the pull-back $F^{(\N - \N')*}_{X^{(\N')}/k} (\mcE)$ is verified to be stable and maximally $F^{(\N')}$-destabilized.
Indeed, suppose that $F^{(\N - \N')*}_{X^{(\N')}/k} (\mcE)$ is not stable. Then, there is a line subbundle $\mcN$ of $F^{(\N - \N')*}_{X^{(\N')}/k} (\mcE)$ such that $\mr{deg}(\mcN) \geq 0$, and we have $F^{(\N')*}_{X/k}(\mcN) \subseteq \mcL$. 
Since $F^{(\N')*}_{X/k}(\mcN)$ is preserved by the connection 
induced by the canonical $\mcD^{(\N' -1)}$-module structure on 
  $F^{(\N)*}_{X/k}(\mcE)$, 
   the composite
\begin{align}
  \mcT _X \otimes F^{(\N')*}_{X/k}(\mcN) \xrightarrow{\mr{inclusion}} \mcT _X \otimes \mcL \xrightarrow{\bar{\phi}} F^{(\N)*}_{X/k} (\mcE) / \mcL
\end{align}
is the  zero map.
This contradicts the fact that $\bar{\phi}$ is an isomorphism.
In particular, the point $a$ lies in $\mr{Dom} (\mr{Ver}^{2}_{\N \Rightarrow 1})$.
\end{rem}

We now  recall the following result,
generalizing ~\cite[Proposition 3.3]{JoXi} (see also ~\cite{JRXY}, ~\cite{JoPa}, \  ~\cite{Wak5}).
Although it was originally stated  in ~\cite[Theorem 36, (ii)]{Wak5} (for the case of {\it pointed} curves) under slightly stronger assumptions,
the current setting allows us to relax them.
In the categorical statements below, a bundle with trivial determinant is understood to be equipped with a chosen trivialization of its determinant, and morphisms are required to preserve this trivialization.

\bpr \label{Prop33}
The following assertions hold.
\begin{itemize}
\item[(i)]
The assignment $\mcE \mapsto \msF^\heartsuit_\mcE$ defined above determines an equivalence of categories
\begin{align} \label{Eq578}
\left(\begin{matrix} \text{the category of} \\ \text{maximally $F^{(\N)}$-destabilized bundles on $X^{(\N)}$} \\ 
\text{of rank $2$ with a trivialized  determinant}
 \end{matrix} \right) \xrightarrow{\sim}
\left(\begin{matrix} \text{the category of dormant} \\ \text{$\mr{SL}_2^{(\N)}$-opers on $X$} \end{matrix} \right),
\end{align}
which can be regarded as a restriction of the equivalence \eqref{Eq577}  under the consideration mentioned in Remark \ref{REm39}.
\item[(ii)]
Let $\varTheta$ be a theta characteristic of $X$.
Then, \eqref{Eq578} restricts to an equivalence of categories
\begin{align} \label{Eq578e}
\left(\begin{matrix} \text{the category of} \\ \text{maximally $F^{(\N)}$-destabilized bundles on $X^{(\N)}$} \\ 
\text{of rank $2$ with a trivialized  determinant} \\
\text{whose pull-back via $F_{X/k}^{(\N)}$ contains} \\
\text{a line subbundle isomorphic to $\varTheta$}
 \end{matrix} \right) \xrightarrow{\sim}
\left(\begin{matrix} \text{the category of dormant} \\ \text{$(\mr{SL}_2^{(\N)}, \varTheta)$-opers on $X$} \end{matrix} \right).
\end{align}
\end{itemize}
Moreover, if $\N'$ is another positive integer with $\N' \leq \N$, then the equivalences 
\eqref{Eq578}  and \eqref{Eq578e}, respectively, commute  with the level-$\N'$ reduction functors  on both sides.
\epr
\begin{proof}
To prove assertion (i),
we   construct an inverse to the assignment  $\mcE \mapsto \msF^\heartsuit_\mcE$.
Let $\msF^\heartsuit := (\mcF, \nabla, \mcL)$ be a dormant $\mr{SL}_2^{(\N)}$-oper on $X$.
Since $\nabla$ has vanishing $p^\N$-curvature,  it follows from the equivalence  \eqref{Eq577} that
$F^{(\N)*}_{X/k}(\mcS ol (\nabla))$ can be identified with $\mcF$ via the natural isomorphism $F^{(\N)*}_{X/k}(\mcS ol (\nabla)) \xrightarrow{\sim} \mcF$.
Note that the vector bundle $\mcS ol (\nabla)$ on $X^{(\N)}$ is stable.
In fact, suppose, on the contrary, that there exists a line subbundle $\mcN$  of $\mcS ol (\nabla)$ of nonnegative degree.
Since $\mr{deg}(F^{(\N)*}_{X/k}(\mcN)) \geq 0 > \mr{deg}(\mcF/\mcL) \left(=1-g \right)$,
the composite 
\begin{align}
F^{(\N)*}_{X/k} (\mcN) \xrightarrow{\mr{inclusion}} F^{(\N)*}_{X/k}(\mcS ol (\nabla)) \xrightarrow{\sim} \mcF  \twoheadrightarrow \mcF/\mcL
\end{align}
 must be the zero map.
This implies  $F^{(\N)*}_{X/k} (\mcN) \subseteq \mcL$, contradicting  the fact that  
$\mcL$ is not closed under $\nabla$ ($= \nabla^{(\N -1)}_{\mcS ol (\nabla), \mr{can}}$ via  the identification  $F^{(\N)*}_{X/k}(\mcS ol (\nabla)) = \mcF$).
Thus, we conclude that $\mcS ol (\nabla)$ is stable.
A similar argument shows that $F^{(\N -1)*}_{X^{(1)}/k} (\mcS ol (\nabla))$ is semistable.
Since $\mr{deg}(\mcL) = g-1$, the vector bundle $\mcS ol (\nabla)$ specifies  an object in  the left-hand side of \eqref{Eq578}.
The resulting  assignment $\msF^\heartsuit \mapsto \mcS ol (\nabla)$ is immediately  verified to be an inverse to  $\mcE \mapsto \msF^\heartsuit_\mcE$.
 This proves assertion (i).

Moreover, assertion (ii) follows from the definition of the assignment $\mcE \mapsto \msF_\mcE^\heartsuit$,  and the remaining portion 
 is a consequence of  the functorial nature of the constructions involved in \eqref{Eq578} and \eqref{Eq578e}.
\end{proof}

\subsection{Comparison of deformation spaces} \label{SS10f}

Let $a$ (resp., $b$) be a $k$-rational point of  $SU_{X^{(\N)}}^2$ (resp., $\mcO p^{^\mr{Zzz...}}_{\N, g}$), and denote by $\mcE$ (resp., $\msF^\spadesuit$) the object classified by  this point.
By fixing a  theta characteristic $\varTheta$ of $X$, we can associate, via  the bijection \eqref{Eq55},  to $\msF^\spadesuit$ a dormant $(\mr{SL}_2^{(\N)}, \varTheta)$-oper $\nabla$ that is normal in the sense of ~\cite[Definition 5.31, (i)]{Wak4} (cf. Remark \ref{Rem44} and ~\cite[Proposition-Definition 5.32]{Wak4}).
Suppose that $\mcE$ corresponds to $(\mcF, \nabla, \varTheta)$, where $\mcF := \mcF_\varTheta$, under   the equivalence of categories \eqref{Eq578}.
 In particular, $\mcE$ is maximally $F^{(\N)}$-destabilized.

We begin with  the following observation, most of which was  already established in the discussion following     ~\cite[Proposition 8.5]{Wak4}.

\ble \label{Lem32}
\begin{itemize}
\item[(i)]
There exists a canonical isomorphism of $k$-vector spaces
\begin{align} \label{Eq230}
\kappa_{\nabla} : H^1 (X^{(\N)}, \mfg_{\mcE}) \xrightarrow{\sim} H^1 (X, \widetilde{\mcT}_{\mcF, B, \nabla}).
\end{align}
\item[(ii)]
Let $\N'$ be another positive 
 integer  with $\N' \leq \N$.
 Denote by  $\nabla'$  the dormant $(\mr{SL}_2^{(\N')}, \varTheta)$-oper $\nabla^{(\N' -1)}$ (i.e., induced from $\nabla$)  and by  $\mcE'$ the maximally $F^{(\N')}$-destabilized bundle $F^{(\N - \N')*}_{X^{(\N')}/k} (\mcE)$ on $X^{(\N')}$ (cf. Remark \ref{Rem33}).
 (Hence, $\nabla'$ corresponds to $\mcE'$ via the equivalence \eqref{Eq578}).
Then,
the following square diagram commutes:
 \begin{align} \label{Ed63df}
\vcenter{\xymatrix@C=46pt@R=36pt{
H^1 (X^{(\N)}, \mfg_{\mcE}) \ar[r]^-{\kappa_\nabla} \ar[d]_-{F^{(\N - \N')*}_{X^{(\N')}/k}} &
H^1 (X, \widetilde{\mcT}_{\mcF, B, \nabla})
 \ar[d]^-{\iota_{\N \Rightarrow \N'}} \\ 
H^1 (X^{(\N')}, \mfg_{\mcE'}) \ar[r]_-{\kappa_{\nabla'}} & 
H^1 (X, \widetilde{\mcT}_{\mcF, B, \nabla'})
 }}
\end{align}
(cf. \eqref{Eq38} for the left-hand vertical arrow and \eqref{Ed63} for the right-hand vertical arrow).
\end{itemize}
\ele
\begin{proof}
We begin with  assertion (i).
Let $\nabla^{(0)}$ denote the injection $\mcT_X \hookrightarrow \widetilde{\mcT}_\mcF$ corresponding to the connection induced by $\nabla$.
 Consider the composite surjection
\begin{align} \label{Eq69}
d'_{\mcF} : \widetilde{\mcT}_\mcF 
\twoheadrightarrow \widetilde{\mcT}_\mcF / \widetilde{\mcT}_{\mcF, B} \xrightarrow{\sim}  \mfg_\mcF /\mfg^0_\mcF \xrightarrow{\sim} \mcT_X
\end{align}
(cf. \eqref{Eq221} for the definitions of $\widetilde{\mcT}_\mcF$ and $\widetilde{\mcT}_{\mcF, B}$).
Here, 
the third arrow is the isomorphism defined in  ~\cite[Eq.\,(617)]{Wak9}, which satisfies $d'_\mcF \circ \nabla^{(0)} = \mr{id}_{\mcT_X}$,  and 
the second arrow
is obtained by applying the snake lemma to 
 the following injection of short exact sequences:
 \begin{align} \label{Eq444}
\vcenter{\xymatrix@C=46pt@R=36pt{
0 \ar[r]  &\mfg_\mcF^0 \ar[r] \ar[d]^-{\mr{inclusion}} & \widetilde{\mcT}_{\mcF, B} \ar[r]^-{d_{\mcF, B}} \ar[d]^-{\mr{inclusion}} & \mcT_X \ar[r] \ar[d]_-{\wr}^-{\mr{id}_{\mcT_X}} & 0
\\
0 \ar[r] & \mfg_\mcF \ar[r] & \widetilde{\mcT}_\mcF \ar[r]_-{d_{\mcF}} & \mcT_X \ar[r] & 0.
  }}
\end{align}
This leads to the definition of  an $\mcO_X$-linear endomorphism 
\begin{align}
\chi := \mr{id}_{\widetilde{\mcT}_\mcF}  - \nabla^{(0)} \circ d_{\mcF} - \nabla^{(0)} \circ d'_{\mcF} : \widetilde{\mcT}_\mcF \rightarrow \widetilde{\mcT}_\mcF
\end{align}
satisfying $d_\mcF \circ \chi = - d'_{\mcF}$
(cf. ~\cite[Lemma  8.4]{Wak4}).
It follows from ~\cite[Proposition 8.5]{Wak4}
that $\chi$ restricts to an automorphism $\chi_{\nabla} : \widetilde{\mcT}_{\mcF, \nabla}\xrightarrow{\sim}\widetilde{\mcT}_{\mcF, \nabla}$ (cf. \eqref{Eq127}), as well as  restricts to an isomorphism 
$\chi_B : \mfg_\mcF \xrightarrow{\sim} \widetilde{\mcT}_{\mcF, B}$.
Hence, 
their intersection $\chi_{B, \nabla} = \chi_{B} \cap \chi_\nabla$ gives rise to 
an isomorphism
\begin{align} \label{EW32}
H^1 (X, \mfg_{\mcF, \nabla})  \xrightarrow{\sim}H^1 (X, \widetilde{\mcT}_{\mcF, B, \nabla}).
\end{align}
Since the equivalence of categories \eqref{Eq577} preserves the formations of ``$\mr{det}$" and ``$\mcH om$", we have $\mfg_{\mcE} \cong \mfg_{\mcF, \nabla}$, which induces   an isomorphism
\begin{align} \label{Eq228}
H^1 (X^{(\N)}, \mfg_\mcE) \xrightarrow{\sim} H^1 (X, \mfg_{\mcF, \nabla}).
\end{align}
Composing  \eqref{EW32} with  \eqref{Eq228}, we obtain an isomorphism $H^1 (X^{(\N)}, \mfg_{\mcE}) \xrightarrow{\sim} H^1 (X, \widetilde{\mcT}_{\mcF, B, \nabla})$, as desired.

Assertion (ii) follows immediately  from  \eqref{Eq10} and the various definitions involved.
Indeed, $\nabla$ and $\nabla'$ induce the same connection $\nabla^{(0)}$, and hence the same morphisms $d'_\mcF$ and $\chi$.
\end{proof}

By composing the isomorphisms  \eqref{Eq36}, \eqref{Eq202}, and \eqref{Eq230},  we conclude  the following assertion.

\bt \label{Prop5}
\begin{itemize}
\item[(i)]
There exists a canonical isomorphism of $k$-vector spaces
\begin{align} \label{Eq29}
\eta_\N : T_a SU^2_{X^{(\N)}} \xrightarrow{\sim} T_b \mcO p_{\N, g}^{^\mr{Zzz...}}.
\end{align}

\item[(ii)]
Let $\N'$ be another positive integer with $\N' \leq  \N$,
and write  $a' := \mr{Ver}^2_{\N \Rightarrow \N'} (a)$, $b':= \Pi_{\N \Rightarrow \N'} (b)$ (hence $a'$ and $b'$ correspond to, respectively,  $\mcE'$ and $\nabla'$ in  the notation of Lemma \ref{Lem32}).
Then,  the following square diagram commutes:
\begin{align} \label{Eq74}
\vcenter{\xymatrix@C=46pt@R=36pt{
 T_a SU^2_{X^{(\N)}}  \ar[r]^-{\eta_\N}_-{\sim} \ar[d]_-{d\mr{Ver}^2_{\N \Rightarrow \N'}|_a} &  T_b \mcO p_{\N, g}^{^\mr{Zzz...}} \ar[d]^-{d\Pi_{\N \Rightarrow \N'}|_b} \\
 T_{a'} SU^2_{X^{(\N')}}  \ar[r]^-{\sim}_-{\eta_{\N'}} &  T_{b'} \mcO p_{\N', g}^{^\mr{Zzz...}}.
 }}
\end{align}
In particular, $\mr{Ver}^2_{\N \Rightarrow \N'}$ is \'{e}tale at $a$ if and only if $\Pi_{\N \Rightarrow \N'}$ is \'{e}tale at $b$.
\end{itemize}
\et

\subsection{Formal completion of  $\mr{Ver}_{\N \Rightarrow \N'}^{2}$} \label{SS88}

We shall examine the structure of  the rational map  $\mr{Ver}_{\N \Rightarrow \N'}^{2}$ restricted over the formal neighborhood of each point classifying  a maximally Frobenius-destabilized bundle.


Let $\mcE$ be  a vector bundle  on a proper, smooth, and connected curve $Y$ over $k$.
A {\bf Higgs field} on $\mcE$ is an $\mcO_Y$-linear morphism $\phi : \mcE \rightarrow \Omega_Y \otimes \mcE$.
We say that  $\mcE$ is {\bf very stable}
if it admits no nonzero nilpotent Higgs field $\phi : \mcE \rightarrow \Omega_Y \otimes \mcE$.
If $Y$ has genus at least two, then every very stable bundle on $Y$ is stable.


\ble \label{LemW1}
For every positive integer $M$, there exists a dense open substack $\mcU_M$ of $\mcM_g$ with the following property:
if $X$ is classified by $\mcU_M$, then, for every integer $\N$ with $1 \leq \N \leq M$ and every dormant $\mr{SL}_2^{(\N)}$-oper $\msF^\heartsuit := (\mcF, \nabla, \mcL)$ on $X$, the vector bundle $\mcS ol (\nabla)$ on $X^{(\N)}$  is very stable.
\ele
\begin{proof}
Choose a totally degenerate stable curve $X_0$ of genus $g$ over $k$, and denote by $x_0$ the $k$-rational point of $\overline{\mcM}_g$ classifying $X_0$.
By ~\cite[Corollary 8.26]{Wak4}, we may find  a common \'{e}tale neighborhood  $S \rightarrow \overline{\mcM}_g$ of $x_0$, with a point $s_0$ above $x_0$, such that all the morphisms $\pi_\N : \mcO p_{\N, g}^{^\mr{Zzz...}} \times_{\overline{\mcM}_g}S\rightarrow S$ for $1 \leq \N \leq M$ obtained by base-change from  $\Pi_\N$ are finite.

Put $S^\mr{sm} := S \times_{\overline{\mcM}_g} \mcM_g$.
For each $\N$,
let $\mcB_\N$ denote the locus of $\mcO p_{\N, g}^{^\mr{Zzz...}} \times_{\overline{\mcM}_g} S^\mr{sm}$ classifying dormant $\mr{PGL}_2^{(\N)}$-opers  represented by  dormant $\mr{SL}_2^{(\N)}$-opers $\msF^\heartsuit := (\mcF, \nabla, \mcL)$ such that $\mcS ol (\nabla)$ admits  a nonzero nilpotent Higgs field.
This condition is independent of the chosen lifting: changing a lifting amounts to tensoring with a dormant line bundle, and the corresponding Cartier descents differ by tensoring with a line bundle.

We verify that $\mathcal B_\N$ is closed. 
\'{E}tale locally on the above stack, the Cartier descents form a family of stable bundles. 
Their trace-free Higgs fields $\phi$ form a vector bundle by Serre duality and the base change theorem.
 The equation
$\mr{det}(\phi)=0$
defines a closed locus in the projectivization of this vector bundle. Its image under the proper projection is precisely $\mcB_\N$,
 since a trace-free rank-$2$ Higgs field is nilpotent if and only if its determinant vanishes. Hence, $\mcB_\N$ is closed.



We claim that the closure $\overline{\mcB}_\N$ in  $\mcO p_{\N, g}^{^\mr{Zzz...}} \times_{\overline{\mcM}_g} S$ does not meet $\pi_\N^{-1} (s_0)$.
Suppose otherwise.
By taking a trait  whose closed point lies in  $\overline{\mcB}_\N \cap \pi_\N^{-1} (s_0)$ and whose generic point lies in $\mcB_\N$,
and  making a finite extension of its fraction field if necessary, we obtain:
\begin{itemize}
\item
a stable curve $X_R$ over a discrete valuation ring $R$, with smooth generic fiber $X_\eta$ and special fiber $X_0$;
\item
a dormant $\mr{GL}_2^{(\N)}$-oper  $\msF^\heartsuit_R := (\mcF_R, \nabla_R, \mcL_R)$ representing the corresponding $\mr{PGL}_2^{(\N)}$-oper on $X_R$ (cf. Remark \ref{Rem9});
\item
a nonzero nilpotent Higgs field  $\phi_\eta$ on the generic fiber $\mcS ol (\nabla_R)_\eta$ of $\mcS ol (\nabla_R)$.
\end{itemize}
Here, the existence of 
such a lifting $\msF^\heartsuit_R$
 follows from ~\cite[Corollary 5.20, (i)]{Wak4}, and the existence of $\phi_\eta$ is unaffected by the choice of this lifting.
 
Let $\Omega_{X_R/R}^\mr{log}$ denote the sheaf of logarithmic $1$-forms on  $X_R/R$ equipped  with the logarithmic structures provided by  ~\cite[Theorem 2.6]{KaFu}.
In particular, $\Omega_{X_R/R}^\mr{log}$ is a line bundle and is isomorphic to the relative dualizing sheaf of $X_R/R$.
Also, we set 
$\mcW_R:=(\Omega_{X_R/R}^\mr{log})^{\otimes p^\N}\otimes\mfg_{\mcF_R}$, 
where $\mfg_{\mcF_R}$ denotes the adjoint bundle of $\mcF_R$, i.e.,  the sheaf of $\mcO_{X_R}$-linear endomorphisms of $\mcF_R$ with vanishing trace.
This vector bundle carries the $\mcD^{(\N -1)}$-module structure $\nabla_{\mcW_R}$ defined as the tensor product of 
the $\mcD^{(\N -1)}$-module structure  on $\mfg_{\mcF_R}$ induced by $\nabla_R$
and 
 the canonical $\mcD^{(\N -1)}$-module structure on $F^{(\N)*}_{X_R/R}(\Omega_{X^{(\N)}_R/R}^\mr{log}) \cong (\Omega_{X_R/R}^\mr{log})^{\otimes p^\N}$
   (cf. \eqref{E445}). 
 


Denote by $\mcW_\eta$ and $\mcW_0$ the generic and special fibers, respectively,  of $\mcW_R$.
Frobenius pullback of $\phi_\eta$ gives a nonzero nilpotent section
$\Phi_\eta \in H^0(X_\eta,\mcW_\eta)$, which is horizontal with respect to $\nabla_{\mcW_R}$.
Its nonvanishing follows from the faithful flatness of the relative Frobenius morphism on the smooth curve $X_\eta$.

Let $t$ be a uniformizer of $R$.
 Properness and flatness of $X_R/R$ imply that
 $H^0(X_R,\mcW_R)$ is a finite torsion-free $R$-module, and
 $H^0(X_R,\mcW_R) \otimes_R \mr{Frac}(R) \cong H^0(X_\eta,\mcW_\eta)$.
 Thus,
  after multiplying $\Phi_\eta$ by a suitable integral power of $t$,  we obtain an extension
$\Phi'_R\in H^0(X_R,\mcW_R)$ that is primitive, i.e., $\Phi'_R\notin tH^0(X_R,\mcW_R)$.
Since the natural morphism 
 $H^0(X_R,\mcW_R)/tH^0(X_R,\mcW_R)
\rightarrow
H^0(X_0,\mcW_0)$ is injective,
the reduction $\Phi'_0$ of $\Phi'_R$ modulo $t$  is nonzero.

Horizontality and nilpotence persist under this extension and reduction:
the corresponding equations hold on the generic fiber and hence on $X_R$,
 since the sheaves in which they take values are $R$-torsion-free. 
 Therefore,
$\Phi'_0
\in
H^0 (X_0, \mcW_0) \setminus \{ 0\}$ is nilpotent and  horizontal  with respect to the special fiber of $\nabla_{\mcW_R}$.

Every component of the normalization of $X_0$, with the preimages of its nodes marked, is a $3$-pointed projective line.
 On each such component, the restriction of $\Phi'_0$ vanishes by 
 ~\cite[Lemma 8.24]{Wak4},
 applied to the restricted  $\mr{GL}_2^{(\N)}$-oper. 
 Consequently, the equality $\Phi'_0=0$ holds, since $X_0$ is reduced and $\mcW_0$ is locally free. This is a contradiction.
 
 Therefore, $\overline{\mcB}_\N\cap\pi_\N^{-1}(s_0)$ is empty.
Since $\pi_\N$ is finite,
 $\pi_\N (\overline{\mcB}_\N)$  is closed in $S$.
  Removing these closed subsets for the various $\N$'s ($1 \leq \N \leq M$)
 gives an open neighborhood of $s_0$.
Its image in $\overline{\mcM}_g$, intersected with $\mcM_g$, is a nonempty dense open substack $\mcU_M$.
It is dense because $\mcM_g$ is irreducible.

By construction, no dormant $\mr{SL}_2^{(\N)}$-oper on a curve classified by $\mcU_M$ with  $1\leq \N\leq M$ has a Cartier descent admitting a nonzero nilpotent Higgs field. 
Thus, $\mcU_M$ satisfies the required condition. 
\end{proof}

\ble \label{LemW2}
Let 
 $R$ be a discrete valuation ring over $k$, 
$Y_R$  a proper, smooth, and geometrically connected curve over $R$,  and $\mcE_R$  a rank-$2$ vector bundle on $Y_R$.
 Denote by $Y_s$ (resp., $Y_\eta$) the fiber of $Y_R$ over the closed point $s$ (resp., the generic point $\eta$) of $\mr{Spec}(R)$, and write $\mcE_s := \mcE_R |_{Y_s}$, $\mcE_\eta := \mcE_R |_{Y_\eta}$.
 Suppose that  $\mcE_s$ is very stable.
Then, every  connection $\nabla_\eta : \mcE_\eta \rightarrow \Omega_{Y_\eta/\eta} \otimes \mcE_\eta$ with vanishing $p$-curvature extends uniquely to a relative connection  
with vanishing $p$-curvature  on $\mcE_R$.
\ele 
\begin{proof}
Let $t$ be a uniformizer of $R$.
By clearing denominators on a finite affine cover of $Y_R$,
there exists a least integer $a\geq 0$ such that
$\nabla'_\eta := t^a \nabla_\eta$
extends to an $R$-linear differential operator of order at most one $\nabla'_R : \mcE_R \rightarrow \Omega_{Y_R/R} \otimes \mcE_R$.
This extension satisfies
$\nabla'_R(fv)
=f\nabla'_R(v)+t^a \cdot df\otimes v$
for local sections $f \in \mcO_{Y_R}$ and $v \in \mcE_R$.

Suppose that $a>0$.
The above Leibniz rule implies that the reduction   of $\nabla'_R$ modulo $t$ defines  an $\mcO_{Y_s}$-linear morphism
$\overline{\nabla}'_s : \mcE_s \rightarrow \Omega_{Y_s/s}\otimes \mcE_s$.
Moreover, 
$\overline{\nabla}'_s \neq 0$
 by the minimality of $a$.
 For a local section $\partial \in \mcT_{Y_R/R}$, 
 denote by $\nabla'_{R, \partial}$ the contraction of $\nabla'_R$ with $\partial$.
Since the $p$-curvature of $\nabla_\eta$ vanishes and $\partial (t)=0$,  the identity 
$({\nabla'_{R, \partial}})^p
=
t^{a(p-1)}\nabla'_{R, \partial^{[p]}}$ holds on the generic fiber.
Both sides are   defined on $Y_R$, so the identity holds on $Y_R$, since $\mcE_R$ is $R$-torsion-free.
 Reducing modulo $t$ yields
$\overline{\nabla}'_{s} (\partial_s)^p = 0$, where $\partial_s$ denotes the restriction of $\partial$ to $Y_s$.
Taking $\partial$ to be a local generator of $\mcT_{Y_R/R}$, we see that $\overline{\nabla}'_s$ is
 a nonzero nilpotent Higgs field on $\mcE_s$, contradicting very stability.
 
 Hence $a=0$, and $\nabla'_R$ is a relative connection extending $\nabla_\eta$.
  Its $p$-curvature vanishes because it vanishes on the generic fiber and $\mcE_R$ is $R$-torsion-free. 
  Uniqueness follows in the same way, by comparing two extensions on the generic fiber.
\end{proof}

Let $\N$ and $\N'$ be positive integers with   $\N' \leq \N$.
By applying Lemmas \ref{LemW1}, \ref{LemW2}, we obtain the following assertion.

\bpr \label{Prop333}
Let $R$ be a discrete  valuation ring over $k$ with fraction field  $K$ and  residue field  $k$. 
 We shall write 
 $\eta$ (resp., $s$)
   for   the generic (resp.,  closed) point of $T := \mr{Spec} (R)$.
Suppose  that $X$ is general in $\mcM_g$.
Consider 
a commutative square diagram of $k$-schemes
\begin{align} \label{Eq1120E}
\vcenter{\xymatrix@C=46pt@R=36pt{
\eta\ar[r]^-{q_\eta} \ar[d]_-{\mr{inclusion}}  & \mr{Dom} (\mr{Ver}_{\N \Rightarrow \N'}^2) \ar[d]^-{\mr{Ver}_{\N \Rightarrow \N'}^2}
\\
T \ar[r]_-{q_T} & SU_{X^{(\N')}}^2
 }}
\end{align}
  such that 
 the vector bundle classified by 
  $q_T (s)$ is   maximally $F^{(\N')}$-destabilized.
Then,
after possibly replacing $R$ with its weakly unramified extension in the sense of
 ~\cite[Tag 0ASG]{Sta}, 
 there exists a unique  morphism $\widetilde{q}_T : T \rightarrow  \mr{Dom} (\mr{Ver}_{\N \Rightarrow \N'}^2)$ making the following diagram commutative:
\begin{align} \label{Eq1120}
\vcenter{\xymatrix@C=46pt@R=36pt{
\eta \ar[r]^-{q_\eta} \ar[d]_-{\mr{inclusion}}  & \mr{Dom} (\mr{Ver}_{\N \Rightarrow \N'}^2) \ar[d]^-{\mr{Ver}_{\N \Rightarrow \N'}^2}
\\
T \ar[r]_-{q_T}  \ar[ur]^{\widetilde{q}_T}& SU_{X^{(\N')}}^2.
 }}
\end{align}
\epr
\begin{proof}
The case $\N=\N'$ is immediate. 
Henceforth, we assume that $\N>\N'$.
It suffices to prove the assertion for   $X \in \mcU_{\N -1}$, where $\mcU_{\N -1}$ is the dense open substack provided by  Lemma \ref{LemW1}.
For each $\Box \in \{ T, \eta, s \}$, we use the notation $(-)_\Box$ to denote base-change to  $\Box$.
For $a \leq b$, write $F_{a, b} : X^{(a)}_\Box \rightarrow X^{(b)}_\Box$ for the $(b-a)$-th  relative Frobenius morphism of $X^{(a)}_\Box$ over $\Box$.
After making the weakly unramified base change allowed in the statement, we may choose stable vector bundles
$\mcE_{\N', R}$ on $X_R^{(\N')}$, $\mcF_\eta$ on $X_\eta^{(\N)}$, both with trivial determinant,  representing  $q_T$ and $q_\eta$, respectively.
The commutativity of the given diagram implies that we may fix an isomorphism $F^*_{\N', \N}(\mcF_\eta) \cong \mcE_{\N', \eta}$.
Put  $\mcG_s := F_{0, \N'}^* (\mcE_{\N', s})$.
Since $\mcE_{\N', s}$ is assumed to  be maximally $F^{(\N')}$-destabilized, 
$F^*_{1, \N'} (\mcE_{\N', s})$ is semistable, and 
$\mcG_s$ has a line subbundle
$\mcL$ of degree  $g-1$.

Consider the  canonical 
$\mcD^{(\N' -1)}$-module structure $\nabla_\mr{can}$ on $\mcG_s$ (cf. \eqref{E445}).
Since $\mcE_{\N', s}$ has trivial determinant, the determinant of $\nabla_{\mr{can}}$ is trivial.
The composite
\begin{align}
\beta : \mcL \hookrightarrow \mcG_s \xrightarrow{\nabla^{(0)}_\mr{can}} \Omega_{X} \otimes \mcG_s \twoheadrightarrow \Omega_X \otimes (\mcG_s/\mcL)
\end{align}
is $\mcO_X$-linear, where the second arrow denotes the connection on $\mcG_s$ induced by $\nabla_\mr{can}$.
We first verify that $\beta$ is an isomorphism.
If $\beta=0$, 
then $\mcL$ is preserved by  $\nabla^{(0)}_\mr{can}$. 
Cartier descent would therefore give a line subbundle $\mcL_1$
of $F^*_{1, \N'} (\mcE_{\N', s})$
 with
$\mcL\cong F_{0,1}^{*}(\mcL_1)$.
Then, we have  $\mr{deg}(\mcL_1) = \frac{1}{p} \cdot \mr{deg}(\mcL) =   \frac{g-1}{p} >0$,  contradicting the semistability of  $F^*_{1, \N'} (\mcE_{\N', s})$.
Consequently, we have  $\beta\neq 0$.
  Since
$\mr{deg} (\mcL)=g-1=\mr{deg} (\Omega_X \otimes (\mcG_s/\mcL))$,
the morphism $\beta$ is an isomorphism.
Thus, $(\mcG_s, \nabla_\mr{can}, \mcL)$, equipped with the induced determinant trivialization,  forms a dormant $\mr{SL}_2^{(\N')}$-oper.

For  each integer $j$ with $\N'\leq j\leq \N$, we define
$\mcE_{j,\eta}:=F_{j, \N}^{*}(\mcF_\eta)$, using the chosen isomorphism when $j=\N'$.
Each $\mcE_{j, \eta}$ has trivial determinant and is stable, since its pull-back along $F_{\N', j}$ is the stable bundle $\mcE_{\N', \eta}$.
We construct their extensions $\mcE_{j,R}$ over $R$ successively, starting with the given $\mcE_{\N',R}$.
We require that each 
$\mcE_{j,R}$ be a family of stable bundles with trivial determinant, and that  
$F_{\N',j}^{*}(\mcE_{j,R}) \cong \mcE_{\N',R}$ compatibly with the prescribed generic identification.
We also require that 
its special fiber be  the Cartier descent of a dormant $\mr{SL}_2^{(j)}$-oper on $X$, whose underlying filtered bundle is $(\mcG_s, \mcL)$ and whose level-$\N'$ reduction is the oper constructed above.
 Suppose that $\mcE_{j,R}$ has been constructed for some $j<\N$.
By Lemma \ref{LemW1}, its special fiber $\mcE_{j,s} \left(:= (\mcE_{j, R})_s \right)$ is very stable.
 The identity
$\mcE_{j,\eta}=F_{j,j+1,\eta}^{*}(\mcE_{j+1,\eta})$
equips $\mcE_{j,\eta}$ with its canonical connection $\nabla_{j, \eta}$ with  vanishing $p$-curvature.
By Lemma \ref{LemW2}, this connection extends uniquely to  a relative connection $\nabla_{j, R}$ over $R$ with vanishing $p$-curvature.
 Relative Cartier descent then gives a vector bundle $\mcE_{j+1,R}$ on $X_R^{(j+1)}$, together with an isomorphism
$F_{j,j+1}^{*}(\mcE_{j+1,R}) \cong
\mcE_{j,R}$,
extending the prescribed generic isomorphism.

Fix a trivialization of $\mr{det}(\mcE_{j, R})$.
With respect to this trivialization, 
 the determinant connection $\mr{det}(\nabla_{j, R})$ is trivial on the generic fiber.
 Indeed,
 $\mr{det}(\mcE_{j +1, \eta})$ is trivial, and any two trivializations of $\mr{det}(\mcE_{j, \eta})$ differ by an element of $K^\times$, whose relative differential is zero.
 Thus, $\mr{det}(\nabla_{j, R})$ agrees with the trivial relative connection over $R$, since the relevant sheaf of homomorphisms is $R$-torsion-free.
 Compatibility of Cartier descent with determinants therefore gives
$\mr{det} (\mcE_{j+1,R})
\cong\mcO_{X_R^{(j+1)}}$.
Furthermore, $\mcE_{j+1,s}$ is stable. Indeed, a line subbundle of nonnegative degree would pull back along $F_{j,j+1,s}$ to a line subbundle of nonnegative degree in the stable bundle $\mcE_{j,s}$.

 Finally, observe that 
$F_{0,j+1}^{*}(\mcE_{j+1,s})\cong F_{0,j}^{*}(\mcE_{j,s}) \cong \mcG_s$.
By transitivity of Cartier descent, the resulting
 $\mcD^{(j)}$-module structure on $\mcG_s$  (cf. \eqref{E445}) restricts to the preceding $\mcD^{(j-1)}$-module  structure,
 and hence to its original $\mcD^{(\N'-1)}$-module structure.
 Together with the determinant trivialization obtained above, this gives a dormant $\mr{SL}^{(j+1)}_2$-oper whose Cartier descent is $\mcE_{j+1,s}$. This completes the induction.
   
  We obtain a family $\mcE_{N,R}$ of stable rank-$2$ bundles with trivial determinant satisfying
$F_{\N',\N}^{*}(\mcE_{\N,R}) \cong \mcE_{\N',R}$ and  $\mcE_{\N, R} |_{X_\eta^{(\N)}}\cong\mcF_\eta$.
It defines a morphism
$\widetilde{q}_T : T\rightarrow SU^{2}_{X^{(\N)}}$.
Since its Frobenius pullback  is the stable family $\mcE_{\N',R}$,
 this morphism factors through
$\mr{Dom}(\mr{Ver}^{2}_{\N\Rightarrow \N'})$
and satisfies  the required property.

Uniqueness follows from separatedness of this domain: two morphisms from the integral scheme $T$ which agree on its generic point agree everywhere.
This completes the proof of this proposition.
\end{proof}


By applying Proposition \ref{Prop333}, we obtain the following assertion.

\bpr \label{Prop19}
Let us take a point $a$ of $SU_{X^{(\N')}}^2$ classifying a maximally $F^{(\N')}$-destabilized bundle.
Denote by  $D$ the spectrum of the completed local ring of $SU_{X^{(\N')}}^2$ at $a$:
\begin{align}
D := \mr{Spec} (\widehat{\mcO}_{SU_{X^{(\N')}}^2, a}) \left(\cong \mr{Spec}(k[\![t_1, \cdots, t_{3g-3}]\!])  \right).
\end{align}
Suppose that $X$ is general in $\mcM_g$.
Then, the  fiber of  $\mr{Ver}^2_{\N \Rightarrow \N'}$ over $D$ decomposes as the disjoint union of  $\mr{deg}(\mr{Ver}^2_{\N \Rightarrow \N'})$ copies of $D$, i.e., 
\begin{align}
\mr{Dom} (\mr{Ver}^2_{\N \Rightarrow \N'}) \times_{SU^2_{X^{(\N')}}} D  \cong   D^{\sqcup \mr{deg}(\mr{Ver}^2_{\N \Rightarrow \N'})}.
\end{align}
In particular, the rational map $\mr{Ver}^2_{\N \Rightarrow \N'}$ is \'{e}tale at all the points classifying maximally $F^{(\N)}$-destabilized bundles.
\epr
\begin{proof}
Denote by $\overline{a}$ the $k$-rational point of $\mcM_g \left(\subseteq \overline{\mcM}_g \right)$  classifying $X$.
Since $X$ is general, it follows from ~\cite[Theorem C, (i)]{Wak4} that 
the projection $\Pi_\N$ is \'{e}tale over $\overline{a}$ and the same  holds  for $\Pi_{\N'}$. 
It follows that the fiber over $\overline{a}$ of  the morphism  $\Pi_{\N \Rightarrow \N'} : \mcO p_{\N, g}^{^\mr{Zzz...}} \rightarrow \mcO p_{\N', g}^{^\mr{Zzz...}}$ is \'{e}tale.
Hence, the second assertion follows from  Theorem \ref{Prop5}, (ii),   and  the comment in  Remark \ref{Rem33}.

Moreover, by  Proposition  \ref{Prop333},
each connected component of  $\mr{Dom} (\mr{Ver}^2_{\N \Rightarrow \N'}) \times_{SU^2_{X^{(\N')}}} D$  contains at least one point mapped to $a$.
By the \'{e}taleness stated in the second assertion, the projection to $D$ is \'{e}tale at this point.
 Since $D$ is strictly henselian, an \'{e}tale neighborhood of this point admits a section over $D$ through the point. 
 The image of this section is open, and it is also closed by separatedness. 
 Hence, the connected component containing this point is isomorphic to $D$.
It follows that $\mr{Dom} (\mr{Ver}^2_{\N \Rightarrow \N'}) \times_{SU^2_{X^{(\N')}}} D$ decomposes as the disjoint union of finitely many copies of $D$, and  the number of its components coincides with $\mr{deg}(\mr{Ver}^2_{\N \Rightarrow \N'})$.
This completes the proof of the first assertion.
\end{proof}

\subsection{Degrees of $\mr{Ver}_{\N \Rightarrow \N'}^{2}$ and $\Pi_{ \N \Rightarrow \N'}$} \label{SS88gge}

One can  apply   Proposition \ref{Prop19} established above   in order to relate  the generic degree of the generalized Verschiebung map  and   the degree of  the projection $\Pi_{\N}$.
As a result, we obtain 
  the following theorem.

\bt \label{Cor30}
Let $\N$ and $\N'$ be positive integers with  $\N' \leq \N$.
Then, the following inequality holds:
\begin{align} \label{Eq234}
\mr{deg}(\Pi_{\N\Rightarrow \N'}) \geq  \mr{deg} (\mr{Ver}^2_{\N \Rightarrow \N'})
 \left(= \mr{deg} (\mr{Ver}^2_{1})^{\N - \N'} \right).
\end{align}
If, moreover, $X$ is general in $\mcM_g$, then this  inequality becomes an equality. 
\et
\begin{proof}
As explained in ~\cite[Section 1]{HoWa},
 the function on $\mcM_g$ given by  $X  \mapsto \mr{deg} (\mr{Ver}_{1}^2)$ (as well as  $X  \mapsto \mr{deg} (\mr{Ver}_{\N \Rightarrow \N'}^2)$) is lower semicontinuous.
Therefore, to obtain   an upper bound of the value $\mr{deg} (\mr{Ver}_{\N \Rightarrow \N'}^2)$, 
we are always free to replace $X$ with a general curve in $\mcM_g$.
Thus, the problem is reduced to proving \eqref{Eq234} with ``$\geq$" replaced by ``$=$" under the assumption that  $X$ is  general in $\mcM_g$.

Choose $X$ sufficiently general so that $\Pi_\N$ and $\Pi_{\N'}$ are finite \'{e}tale over a neighborhood of the point $q$ classifying $X$, and Proposition \ref{Prop19} applies. 
Let $b'$ be any point of $\mcO p^{^\mr{Zzz...}}_{\N',g}$ over $q$.
 Fix a theta characteristic $\varTheta$, and let $a\in SU^{2}_{X^{(\N')}}$ correspond to $b'$ under \eqref{Eq55} and \eqref{Eq578e}.
Compatibility of these correspondences with level reduction identifies the geometric fiber of $\Pi_{\N\Rightarrow \N'}$ over $b'$ with that of $\mr{Ver}^{2}_{\N\Rightarrow \N'}$ over $a$. By Proposition \ref{Prop19}, both fibers consist of
$\mr{deg}(\mr{Ver}^{2}_{\N\Rightarrow \N'})$ points.
This number is independent of $b'$. Hence, over a suitable dense open substack of $\mcM_g$, 
the morphism $\Pi_{\N\Rightarrow \N'}$ has this constant degree. This proves the desired equality.
\end{proof}

\bco[cf. Theorem \ref{ThA}] \label{Cor88}
For a positive integer $\N$, the following inequality holds:
\begin{align}  \label{Eq245}
\mr{deg}(\Pi_{\N}) \geq   \frac{\mr{deg}(\mr{Ver}_{1}^{2})^{\N -1} \cdot p^{g-1}}{2^{2g-1}} \cdot \sum_{\theta =1}^{p-1} \frac{1}{\sin^{2g-2} \big(\frac{\pi \cdot \theta}{p} \big)}.
\end{align}
If, moreover, $X$ is general in $\mcM_g$,
then this  inequality becomes an equality.
\eco
\begin{proof}
As discussed in the proof of Theorem \ref{Cor30},
the problem is reduced to proving  \eqref{Eq245} when ``$\geq$" is replaced by ``$=$"  under the assumption that $X$ is general in $\mcM_g$.
Let us  observe that 
\begin{align} 
\mr{deg}(\Pi_{\N}) &= \mr{deg}(\Pi_{1} \circ \Pi_{2 \Rightarrow 1} \circ \Pi_{3 \Rightarrow 2} \circ \cdots \circ \Pi_{\N \Rightarrow \N -1})  \\
&= \mr{deg}(\Pi_1) \cdot \prod_{i=2}^{\N} \mr{deg}(\Pi_{i \Rightarrow i-1}) \notag \\
& = \mr{deg}(\Pi_1) \cdot \prod_{i=2}^{\N}\mr{deg}(\mr{Ver}^2_{i \Rightarrow i-1})  \notag \\
& = \mr{deg}(\Pi_1) \cdot \prod_{i=2}^{\N}\mr{deg}(\mr{Ver}^2_{1})  \notag \\
&= \mr{deg}(\Pi_1) \cdot   \mr{deg}(\mr{Ver}^2_{1})^{\N-1} \notag \\
& =  \left(\frac{p^{g-1}}{2^{2g-1}} \cdot \sum_{\theta =1}^{p-1} \frac{1}{\sin^{2g-2} \big(\frac{\pi \cdot \theta}{p} \big)}\right) \cdot \mr{deg}(\mr{Ver}_{1}^{2})^{\N -1},
\end{align}
where the third equality follows from Theorem \ref{Cor30} and  the last equality follows from ~\cite[Theorem A]{Wak1}.
Here, the formula for $\mr{deg}(\Pi_1)$ extends to every odd prime $p$.
 Indeed, both sides are polynomial in $p$, by  ~\cite[Theorem 2.1]{LiOs} and 
 ~\cite[Section 6.2, (2)]{Wak1},
 respectively. 
 Since  ~\cite[Theorem A]{Wak1}
 identifies these polynomials for all primes $p>2g-2$, they agree identically.
This completes the proof of this assertion.
\end{proof}

\begin{rem} \label{Rem56}
In ~\cite[Conjecture 10.29]{Wak4}, the second author  proposed  a conjecture predicting an explicit  connection between the generic degree of the generalized Verschiebung map on the moduli  space of rank $n (>1)$ stable bundles and the degree of the moduli stack of dormant  $\mr{PGL}_n^{(\N)}$-opers.
The above theorem provides 
 an affirmative answer to this conjecture  in the case $n= 2$; indeed, 
\eqref{Eq245} yields  
 the following equality for a general $X$:
\begin{align}
\mr{deg}(\mr{Ver}_1^2) = \lim_{\N \to \infty} \mr{deg}(\Pi_\N)^{1/\N}.
\end{align}
We expect that the general case of this conjecture might  be resolved  by arguments similar to those developed in this paper after proving  the generic \'{e}taleness of the moduli stack of dormant $\mr{PGL}_n^{(\N)}$-opers.
\end{rem}

\vspace{10mm}
\section{Computing  the generic degree of the generalized Verschiebung} \label{Sgr2}

In this final section, 
we recall a combinatorial description of
dormant $\mr{PGL}_2^{(\N)}$-opers on totally degenerate curves established  in ~\cite[Section 10]{Wak4}.
Moreover, we  apply this previous work together with  the results obtained so far to
provide a procedure for explicitly computing the generic degree of the generalized Verschiebung.
For the related discussions in the case $\N =1$, we refer the reader to ~\cite{LiOs}, ~\cite{Moc}, ~\cite{Wak2}, ~\cite{Wak9}, and ~\cite{Wak6}.

\subsection{Balanced $(p, \N)$-edge numberings} \label{SS5}

For each positive integer $\N$, 
we shall set 
\begin{align} \label{Eq260}
{^\dagger}C_\N
\end{align}
(cf. ~\cite[Section 10.2]{Wak4}) to be the set of  triples of nonnegative integers $(s_1, s_2, s_3)$ satisfying the following conditions:
\begin{itemize}
\item
$\sum_{i=1}^3 s_i \leq p^\N -2$ and $|s_2 -s_3| \leq s_1 \leq s_2 + s_3$;
\item
For every positive integer $\N' < \N$, there exists a triple  $(s'_1, s'_2, s'_3)$ with $s'_i \in \{ [s_i]_{\N'}, p^{\N'} -1 - [s_i]_{\N'} \}$ for $i=1,2,3$ such  that  $\sum_{i=1}^3 s'_i \leq p^{\N'} -2$ and $|s'_2 - s'_3| \leq s'_1 \leq s'_2 + s'_3$.
\end{itemize}
Here, for each nonnegative integer $a$, we set $[a]_{\N'}$
to be the remainder obtained by dividing $a$ by $p^{\N'}$.
Note that the second condition can be equivalently expressed as requiring that, for every positive  integer $\N' < \N$,
the triple $(\hat{s}_1, \hat{s}_2, \hat{s}_3)$ satisfies 
 $\sum_{i=1}^3 \hat{s}_i \leq p^{\N'} -2$ and $|\hat{s}_2 - \hat{s}_3| \leq \hat{s}_1 \leq \hat{s}_2 + \hat{s}_3$, where
 \begin{align}
 \hat{a} := \frac{p^{\N'}-1}{2} - \left| a - p^{\N'} \cdot \left\lfloor \frac{a}{p^{\N'}}\right\rfloor  - \frac{p^{\N'}-1}{2} \right|
 \end{align}
  for each nonnegative integer $a$.

\begin{rem} \label{Rem100}
In ~\cite{Wak4} (and ~\cite{Moc} for $\N =1$), a characterization   
 of dormant $\mr{PGL}_2^{(\N)}$-opers on a $3$-pointed projective line was provided in terms of specific combinations  of nonnegative integers.
This characterization  is derived   by examining  the radii, which are encoded by the local exponents modulo $p^\N$, of such opers at the marked points.
As a consequence of the theory of diagonal reduction/lifting and Gauss hypergeometric differential operators in characteristic $p^\N$,  
we obtain  a bijective correspondence between dormant  $\mr{PGL}_2^{(\N)}$-opers on this curve and elements of ${^\dagger}C_\N$   (cf. ~\cite[Theorem 10.16]{Wak4}).
\end{rem}

Next, let us fix 
a finite connected trivalent graph $\mbG$ (without open edges), in the sense of ~\cite[Definition 6.5, (i)-(iii)]{Wak4};
this  consists of a vertex set $V$ and  an edge set  $E$, each of whose elements has cardinality $2$ (with $e \cap e' = \emptyset$ if  $e \neq e'$),  and a map $\zeta : \coprod_{e \in E} e \rightarrow V$.
In particular,  we have  $\sharp (V), \sharp (E) < \infty$.
We refer to any element $b$ of $e \in E$ as a branch of  $e$.
For each vertex  $v \in V$, we shall write $B_v := \zeta^{-1} (v)$, which is  the set of branches 
   incident to $v$ (hence $\sharp (B_v) = 3$).
Assume  further that $\mbG$ has genus $g$, i.e., is of type $(g, 0)$, in the sense of   ~\cite[Definition 6.5, (iv)]{Wak4}.
That is, when we regard $\mbG$ as a topological space in the usual way (cf. ~\cite[Remark 7.3]{Wak9}), this assumption requires that  the first Betti number $b_1 (\mbG) := \mr{dim}_\mbQ (H_1 (\mbG, \mbQ))$ coincides with $g$.

\bde \label{Def399}
A {\bf balanced $(p, \N)$-edge numbering} on $\mbG$ is a collection
$(a_{e})_{e \in E}$
of nonnegative integers  indexed by the elements of $E$ such that, for each vertex $v \in V$, the triple  $(a_{e_b})_{b \in B_v}$  belongs to ${^\dagger}C_\N$, where $e_b$ denotes the unique edge satisfying $b \in e_b$. 
\ede

We denote the set of all balanced $(p, \N)$-edge numberings  on $\mbG$ by 
\begin{align} \label{Eq111}
\mr{Ed}_{p, \N, \mbG}.
\end{align}
This set contains  finitely many elements, and its cardinality can be computed by explicitly  enumerating balanced $(p, \N)$-edge numberings.
For example, when $\mbG$ has genus $3$, we have 
\begin{align}
&\sharp (\mr{Ed}_{3, 2, \mbG}) = 49,  \ \ \ 
\sharp (\mr{Ed}_{5, 2, \mbG}) = 11775,  \ \ \ 
\sharp (\mr{Ed}_{7, 2, \mbG}) = 542626,  \ \ \ 
\sharp (\mr{Ed}_{11, 2, \mbG}) = 107098915, \\
& \sharp (\mr{Ed}_{3, 3, \mbG}) = 2401, \ \ \ 
\sharp (\mr{Ed}_{5, 3, \mbG}) = 9243375, \ \ \ 
\sharp (\mr{Ed}_{7, 3, \mbG}) = 3004520162, \\
& \sharp (\mr{Ed}_{3, 4, \mbG}) = 117649, \ \ \ 
\sharp (\mr{Ed}_{5, 4, \mbG}) = 7256049375.
\end{align}

(For details of the following discussion, we refer to ~\cite[Section 10.5]{Wak4}.)
Recall that a stable curve is said to be {\bf totally degenerate} if it is obtained by gluing together finitely many copies of a $3$-pointed projective line, along their marked points to form nodal points (cf. ~\cite[Definition 6.26]{Wak4} for the precise definition).
One can find a totally degenerate curve  $X_\mbG$ of genus $g$ whose dual graph is given by  $\mbG$.
In particular, the vertices in $\mbG$ represent  irreducible components, and the edges represent  nodal points.

Since the projection $\Pi_{\N}$  is finite and \'{e}tale over an \'{e}tale neighborhood  of the point classifying $X_\mbG$
   (cf. ~\cite[Theorem C, (i)]{Wak4}),
its degree  $\mr{deg} (\Pi_\N)$ coincides with the total number of dormant $\mr{PGL}_2^{(\N)}$-opers on $X_\mbG$.
Moreover, each such dormant $\mr{PGL}_2^{(\N)}$-oper  can be constructed  by gluing dormant $\mr{PGL}_2^{(\N)}$-opers on the components of its normalization, each regarded as a $3$-pointed projective line, 
in accordance with the attachment of underlying curves
at the marked points.
As noted  in Remark \ref{Rem100},
these componentwise opers  are classified by their radii, which are identified with  elements of ${^\dagger}C_\N$.
Therefore,  a global  dormant $\mr{PGL}_2^{(\N)}$-oper on $X_\mbG$ corresponds to an edge numbering such that, for every vertex, the numbers on incident edges form a triple in 
  ${^\dagger}C_\N$.
This gives  a bijective correspondence  between dormant $\mr{PGL}_2^{(\N)}$-opers on $X_\mbG$ and elements of $\mr{Ed}_{p, \N, \mbG}$.
Consequently, we obtain   the following equality:
\begin{align} \label{Eq110}
\mr{deg}(\Pi_{\N}) = \sharp  (\mr{Ed}_{p, \N, \mbG})
\end{align}
(cf. ~\cite[Theorem E]{Wak4}).
This fact
and Corollary \ref{Cor88}
 together imply   the following assertion.

\bt[cf. Theorem \ref{ThE}, (i)] \label{Th3} 
Let us keep the above notation.
Suppose that  
 $X$ is general in $\mcM_g$.
Then,
the following equality holds:
\begin{align} \label{Eq288}
\mr{deg}(\mr{Ver}^2_1) = \frac{\sharp (\mr{Ed}_{p, \N+1, \mbG})}{\sharp (\mr{Ed}_{p, \N, \mbG})}
 \left(= \frac{\sharp (\mr{Ed}_{p, 2, \mbG})\cdot 2^{2g-1}}{p^{g-1}}  \cdot \left(\sum_{\theta =1}^{p-1} \frac{1}{\sin^{2g-2} \big(\frac{\pi \cdot \theta}{p} \big)}\right)^{-1} \right). 
\end{align}
\et 

\subsection{Explicit computations for $3\leq g\leq 10$} \label{SS588}

In this final subsection, we consider a polynomial expression for $\mr{deg}(\mr{Ver}_1^2)$.
According to  ~\cite[Proposition 10.26]{Wak4},  the set  $\mr{Ed}_{p, \N, \mbG}$ admits a description in terms of lattice points in dilations of  certain rational constructible sets.
 In particular, one may apply Ehrhart's  theory of lattice-point enumeration for rational polytopes to this setting.
As shown  in  ~\cite[Theorem  10.27]{Wak4},
 there exists a quasi-polynomial $H_{\N} (t)$ of degree $(3g-3) \cdot \N$ with coefficients in $\mbQ$ whose   minimum period divides $4$ and which satisfies  
 \begin{align}
 \sharp (\mr{Ed}_{p, \N, \mbG}) = H_{\N} (p)
 \end{align}
  for every odd prime $p$.

Since $\sharp (\mr{Ed}_{p, 1, \mbG}) \left(= H_{1} (p) \right)$ has been already computed in  prior work,
it follows from Theorem \ref{Th3}  that determining  $\mr{deg} (\mr{Ver}_1^2)$ reduces to 
 identifying  the coefficients of $H_{2} (t)$.
 To this end, 
 it suffices to find  the number of balanced $(p, 2)$-edge numberings on $\mbG$ for some values of $p$, and then solve the resulting system of linear equations.
 
 The following two assertions may slightly reduce the computational burden.

\bpr \label{Th54}
 The leading coefficient of $H_{1} (t)$ (resp.,  $H_{2} (t)$) coincides with 
 $\frac{(-1)^g \cdot   B_{2g-2}}{2 \cdot (2g-2)!}$ (resp., $\frac{2^{3g-5} \cdot  B^2_{2g-2}}{((2g-2)!)^2}$),
 where $B_{2g-2}$ denotes the $(2g-2)$-th Bernoulli number, i.e., $B_2 = \frac{1}{6}$, $B_4 = -\frac{1}{30}$, and $B_6 =  \frac{1}{42}$, etc.
 \epr
 \begin{proof}
 Let us consider the rational polytope $\mcP_0$ in  $\mbR^3$ consisting of all points $(x, y, z)$ satisfying the conditions
\begin{align}
0 \leq x,  \ \ \ 0 \leq y,  \ \ \ 0 \leq z, \ \ \ x + y + z \leq 1,  \ \  \text{and} \ \ \   |y- z| \leq x \leq y + z.
\end{align} 
In particular, $\mcP_0$ is contained in $\left[0, \frac{1}{2} \right]^3 \left(\subseteq \mbR^3 \right)$.
Denote by 
$\mbR^E$ the space of all real-valued functions $E \rightarrow \mbR$ on the edge set $E$; it forms a $(3g-3)$-dimensional $\mbR$-vector space, and the subset of integer-valued functions $\mbZ^E$ forms its lattice.
Each element of $\mbR^E$ (resp., $\mbZ^E$) may be regarded as  a collection $(s_i)_{i \in E}$ with $s_i \in \mbR$ (resp., $s_i \in \mbZ$).

We now define 
$\mcP$ to be  the subset of $\mbR^E$ consisting of collections  $(s_i)_{i \in E}$
such that, for each $v \in V$,
the triple $(s_{e_b})_{b \in B_v}$ lies in $\mcP_0$.
As established  in ~\cite[Introduction, Theorem 1.3]{Moc} (or  ~\cite{LiOs}, ~\cite{Wak2}), 
there exists a natural  bijection  
\begin{align}
\mr{Ed}_{p, 1, \mbG} \xrightarrow{\sim} (p-2) \mcP \cap \mbZ^E,
\end{align}
where  for a subset $\mcQ$ in $\mbR^E$ and a nonnegative integer $m$ we denote by $m \mcQ$ the dilation of $\mcQ$ by the factor $m$.
Thus, $H_1 (t)$ coincides,  up to a shift of  the variable by $2$,   with the Ehrhart quasi-polynomial associated to $\mcP$. 
By the discussion in ~\cite[Section 6.2, (2)]{Wak1},
 their  leading coefficients are given by   
$\frac{(-1)^g \cdot B_{2g-2}}{2 \cdot (2g-2)!}$.
This  completes the proof of the non-resp'd portion.

Next, 
for each $s \in \mbR$, we set $\langle s \rangle_+ := s$ and $\langle s \rangle_- := 1- s$.
We also denote by $\mr{Sgn}_{2, \mbG}$
the set of collections $\mfa := (\mfa_i)_{i \in E}$ such that $\mfa_{i} \in \{ 1, -1 \}$ for every $i \in E$.
Each such collection $\mfa$ determines the  bijective transformation $\delta_\mfa : \mbR^E \xrightarrow{\sim} \mbR^E$
given by $(s_i)_{i} \mapsto (\langle s_i \rangle_{\mfa_i})_i$. 
According to ~\cite[Proposition 10.26, (ii)]{Wak4}, 
there exist 
   subsets
$\mcP_{1} (\mfa)$  and $\mcP_{2} (\mfa)$ (for $\mfa \in \mr{Sgn}_{2, \mbG}$)
of  $\mbR^E$ satisfying the following properties:
\begin{itemize}
\item[(a)]
For  every $\mfa \in \mr{Sgn}_{2, \mbG}$,
the closure $\overline{\mcP_{1} (\mfa)}$ (resp., $\overline{\mcP_{2} (\mfa)}$) of $\mcP_{1} (\mfa)$ (resp.,  $\mcP_{2} (\mfa)$) is a  convex rational polytope satisfying 
\begin{align}
\overline{\mcP_{1} (\mfa)} \supseteq \mcP_{1} (\mfa) \supseteq \overline{\mcP_{1} (\mfa)} \setminus \partial \overline{\mcP_{1} (\mfa)},
\ \left(\text{resp.,} \ 
\overline{\mcP_{2} (\mfa)} \supseteq \mcP_{2} (\mfa) \supseteq \overline{\mcP_{2} (\mfa)} \setminus \partial \overline{\mcP_{2} (\mfa)} \right),
\end{align}
where $\partial \overline{\mcP_{1} (\mfa)}$ (resp., $\partial \overline{\mcP_{2} (\mfa)}$) denotes the boundary of  $\overline{\mcP_{1} (\mfa)}$ (resp., 
$\overline{\mcP_{2} (\mfa)}$);
\item[(b)]
For every $\mfa \in \mr{Sgn}_{2, \mbG}$,
we have $\delta_\mfa (\overline{\mcP_{1} (\mfa)}) = \overline{\mcP_{2} (\mfa)} = \mcP$.
\item[(c)]
There exists a bijective correspondence
\begin{align} \label{Eq409}
\mr{Ed}_{p, 2, \mbG} \cong \coprod_{\mfa \in \mr{Sgn}_{2, \mbG}} \left((p-1) \mcP_{1} (\mfa)\cap \mbZ^{E} \right) \times \left( p \mcP_{2} (\mfa) \cap \mbZ^E\right).
\end{align}
\end{itemize}

Note that   $\mcR :=  \coprod_{\mfa \in \mr{Sgn}_{2, \mbG}} \overline{\mcP_{1} (\mfa)} \times \overline{\mcP_{2} (\mfa)}$ is a finite disjoint union of rational polytopes, each lying in $\mbR^{E \sqcup E}$.
We define its Ehrhart quasi-polynomial $H_\mcR (t)$ to be the sum of the Ehrhart quasi-polynomials of these polytopes.
It then follows from the non-resp'd assertion and the property (b) above  that  the leading coefficient  of 
$H_\mcR (t)$
 coincides with 
\begin{align} \label{Eq223}
\sharp (\mr{Sgn}_{2, \mbG}) \cdot  \left(\frac{(-1)^g \cdot B_{2g-2}}{2 \cdot (2g-2)!} \right)^2 =  \frac{2^{3g-5} \cdot B_{2g-2}^{2}}{((2g-2)!)^2}.
\end{align}
On the other hand,  by the properties (a) and (c) described above,
the leading coefficient of $H_2 (t)$  is verified to be the same as that of $H_\mcR (t)$.
This finishes the proof of the resp'd portion.
 \end{proof}

 Applying the above proposition, we obtain the following result.

\bt[cf. Theorem \ref{ThE}, (ii)] \label{Th56}
Suppose that $X$ is general in $\mcM_g$.
Then, there exists a quasi-polynomial $Q (t)$ with coefficients in $\mbQ$, of degree $3g-3$ and period dividing $4$,  such that 
the generic degree $\mr{deg}(\mr{Ver}_1^2)$ (in characteristic $p$) is given by the value  $Q (p)$ for every odd prime $p$.
Moreover, the leading coefficient of $Q (t)$ coincides with
   $\frac{(-1)^g \cdot  2^{3g-4} \cdot B_{2g-2}}{(2g-2)!}$.
 \et
\begin{proof}
For $r\in\{1,3\}$,
let $A_r(t)$ and $B_r(t)$ denote the constituents of $H_1(t)$ and $H_2(t)$, respectively, on the residue class $r$ modulo $4$. Polynomial division gives
$B_r(t)=A_r(t)Q_r(t)+R_r(t)$ for some 
$Q_r (t),R_r(t)\in\mbQ [t]$ with
$\mr{deg} (R_r (t))<3g-3$.

Choose a positive integer $c_r$ such that $c_rQ_r(t)\in\mbZ [t]$.
 For every odd prime $p\equiv r$ (mod $4$),
  Theorem \ref{Th3} yields
$\frac{c_rR_r(p)}{A_r(p)}=c_r \mr{deg}(\mr{Ver}_1^2)-c_r Q_r(p)\in\mbZ$.
Since $\mr{deg} ( R_r (t))<\mr{deg} ( A_r (t))$, the left-hand side tends to zero as $p$ tends to infinity through these primes. 
It therefore vanishes for all sufficiently large such primes. 
There are infinitely many primes in each of the residue classes $1$ and $3$ modulo $4$, so $R_r (t)=0$.

Thus, $B_r (t)=A_r (t)Q_r (t)$ for $r=1,3$.
Define a quasi-polynomial $Q (t)$ of period dividing $4$ by taking $Q_1(t)$, $Q_3(t)$ as its odd constituents and, for example, $Q_1(t)$ as both of its even constituents.
 Then Theorem \ref{Th3} gives
$Q(p)=\mr{deg}(\mr{Ver}_1^2)$
for every odd prime $p$.
 By Proposition \ref{Th54}, each $Q_r (t)$ has degree $3g-3$,
  and its leading coefficient is
\begin{align}
\frac{
\frac{2^{3g-5} \cdot B_{2g-2}^{2}}{((2g-2)!)^2}
}{
\frac{(-1)^g \cdot B_{2g-2}}{2 \cdot (2g-2)!}
}
=
\frac{(-1)^g\cdot 2^{3g-4} \cdot B_{2g-2}}{(2g-2)!}.
\end{align}
This proves the assertion.
\end{proof}

Let us now focus on the cases $3 \leq g \leq 10$.
To begin with, we observe that the notion of a balanced $(p, \N)$-edge numbering can also  be defined  even when $p$ is replaced by   any odd integer $q \geq 3$.
The combinatorial construction underlying \eqref{Eq409}, based on base-$q$ expansions and the defining inequalities, applies unchanged to these values.
Consequently, the Ehrhart quasi-polynomial $H_2 (t)$ satisfies $H_2 (q) = \sharp (\mr{Ed}_{q, 2, \mbG})$ for every odd integer $q \geq 3$, 
so the coefficients of its odd constituents can be determined using sufficiently many small odd integers, not necessarily primes.

For instance, when $g=3$, by  explicitly counting balanced $(p, 2)$-edge numberings for a genus-$3$ trivalent graph  $\mbG$,
we obtain  
\begin{align}
& \sharp (\mr{Ed}_{3, 2, \mbG})  = 49, \ \ \ 
\sharp (\mr{Ed}_{5, 2, \mbG})  =11775, \ \ \ 
\sharp (\mr{Ed}_{7, 2, \mbG})  = 542626, \\ 
& \sharp (\mr{Ed}_{9, 2, \mbG})  =10108638, \ \ \  
\sharp (\mr{Ed}_{11, 2, \mbG})  = 107098915, \ \ \
\sharp (\mr{Ed}_{13, 2, \mbG}) = 773117709, \\
& \sharp (\mr{Ed}_{15, 2, \mbG})  = 4229656900, \ \  \ 
\sharp (\mr{Ed}_{17, 2, \mbG})  = 18767168700,  \ \ \ 
\sharp (\mr{Ed}_{19, 2, \mbG})  =70695102549,  \\
& \sharp (\mr{Ed}_{21, 2, \mbG})  = 233505804763, \  \ \ 
\sharp (\mr{Ed}_{23, 2, \mbG})  = 692440249446,  \ \ \  
\sharp (\mr{Ed}_{25, 2, \mbG})  = 1876599156250,
\end{align}
and so forth.
The two odd constituents of $Q(t)$ have degree $6$ and leading coefficient $\frac{2}{45}$ (cf.  Theorem \ref{Th56}).
 The six values in each of the residue classes $1$ and $3$ modulo $4$ therefore determine these constituents uniquely. 
 Interpolation gives the same polynomial in both cases:
\begin{align}
Q_1(t)=Q_3(t)=\frac{t^2(2t^4+5t^2+38)}{45}.
\end{align}
These computations  yield  the following result. The remaining cases $4\leq g\leq 10$ can be worked out by the same method; see Appendix  \ref{AppPoly}.

\bt[cf. Theorem \ref{ThG}] \label{Th446}
Suppose that  
 $X$ is general in $\mcM_g$.
Then, for each integer $g$ with $3 \leq g \leq 10$, the generic degree $\mr{deg}(\mr{Ver}^2_1)$ is given by a polynomial in $p$.
See Appendix \ref{AppPoly}
  for the explicit list.
\et


\vspace{10mm}
\appendix
\section{Polynomial formulas for $\mr{deg}(\mr{Ver}^2_1)$ for $3\leq g\leq 10$} \label{AppPoly}
Suppose that $X$ is general in $\mcM_g$.
In this appendix, we list explicit polynomial expressions in $p$ for $\mr{deg}(\mr{Ver}^2_1)$ in the case $3\leq g\leq 10$, obtained by counting balanced $(p,2)$-edge numberings.
\begin{align*}
\intertext{$g=3$:}
\mr{deg}(\mr{Ver}^2_1) &= \frac{p^2(2p^4 + 5p^2 + 38)}{45}. \\
\intertext{$g=4$:}
\mr{deg}(\mr{Ver}^2_1) &= \frac{p^3(8p^6 + 21p^4 + 42p^2 + 874)}{945}. \\
\intertext{$g=5$:}
\mr{deg}(\mr{Ver}^2_1) &= \frac{p^4(24p^8 + 80p^6 + 147p^4 + 270p^2 + 13654)}{14175}. \\
\intertext{$g=6$:}
\mr{deg}(\mr{Ver}^2_1) &= \frac{p^5(32p^{10}+132p^8+286p^6+451p^4+792p^2+91862)}{93555}. \\
\intertext{$g=7$:}
\mr{deg}(\mr{Ver}^2_1) &= \frac{p^6}{638512875}\Big(44224p^{12} + 218400p^{10} + 558558p^8 + 993850p^6 + 1438437p^4 \\
&\hspace{18mm}+2457000p^2 + 632802406\Big). \\
\intertext{$g=8$:}
\mr{deg}(\mr{Ver}^2_1) &= \frac{p^7}{273648375}\Big(3840p^{14} + 22112p^{12} + 65520p^{10} + 133419p^8 + 211640p^6 \\
&\hspace{18mm}+ 289926p^4 + 486000p^2 + 272435918\Big). \\
\intertext{$g=9$:}
\mr{deg}(\mr{Ver}^2_1) &= \frac{p^8}{488462349375}\Big(1388928p^{16} + 9139200p^{14} + 30824128p^{12} + 71073600p^{10} \\
&\hspace{18mm}+ 126249123p^8 + 185684200p^6 + 244813158p^4 + 404838000p^2 + 487388339038\Big). \\
\intertext{$g=10$:}
\mr{deg}(\mr{Ver}^2_1) &= \frac{p^9}{194896477400625}\Big(112299520p^{18} + 831273408p^{16} + 3145141440p^{14} + 8109940848p^{12} \\
&\hspace{18mm}+ 16038878310p^{10} + 25993141083p^8 + 36233012730p^6 + 46442359908p^4 \\
&\hspace{18mm}+ 76014288000p^2 + 194683557065378\Big).
\end{align*}

\subsection*{Acknowledgements}
We are grateful to the referee for a careful reading of the manuscript and for many valuable comments and suggestions, which helped us improve both the arguments and the exposition.
The first author is a JSPS Research Fellow and was supported by the Grant-in-Aid for JSPS Fellows (JSPS KAKENHI Grant Number JP24KJ1651).
 The second author was partially supported by 
 JSPS KAKENHI Grant Number 25K06933.
Also, 
this paper is based on a presentation given by the second author during a lecture series entitled ``Oper, Dormant Oper, $p$-adic Teichm\"{u}ller theory," held at Kyushu University in March 2025. We would like to express our sincere gratitude to Professor Atsushi Katsuda, the organizer of the series, for providing us with this opportunity.

\end{document}